\documentclass[11pt]{article}
\usepackage{epsf}
\usepackage{amssymb}
\usepackage{theorem}
\usepackage{latexsym}
\font\tenmsb=msbm10
\font\sevenmsb=msbm7
\newfam\msbfam 
\textfont\msbfam=\tenmsb
\scriptfont\msbfam=\sevenmsb

\def\Bbb#1{\fam\msbfam\relax#1}

\newtheorem{thm}{Theorem}[section]
\newtheorem{prop}[thm]{Proposition}
\newtheorem{cor}[thm]{Corollary}
\newtheorem{lem}[thm]{Lemma}
\newtheorem{conj}[thm]{Conjecture}
\newtheorem{exa}[thm]{Example}
\newtheorem{defn}[thm]{Definition}
\newtheorem{rem}[thm]{Remark}

\newcommand{\ben}{\begin{enumerate}}
\newcommand{\een}{\end{enumerate}}
\newcommand{\ble}{\begin{lem}}
\newcommand{\ele}{\end{lem}}

\newcommand{\bre}{\begin{rem}}
\newcommand{\ere}{\end{rem}}

\newcommand{\bthm}{\begin{thm}}
\newcommand{\ethm}{\end{thm}}
\newcommand{\bpr}{\begin{prop}}
\newcommand{\epr}{\end{prop}}
\newcommand{\bco}{\begin{cor}}
\newcommand{\eco}{\end{cor}}
\newcommand{\bcon}{\begin{conj}}
\newcommand{\econ}{\end{conj}}
\newcommand{\bde}{\begin{defn}}
\newcommand{\ede}{\end{defn}}
\newcommand{\bex}{\begin{exa}}
\newcommand{\eex}{\end{exa}}
\newcommand{\barr}{\begin{array}}
\newcommand{\earr}{\end{array}}
\newcommand{\btab}{\begin{tabular}}
\newcommand{\etab}{\end{tabular}}
\newcommand{\beq}{\begin{equation}}
\newcommand{\eeq}{\end{equation}}
\newcommand{\bea}{\begin{eqnarray*}}
\newcommand{\eea}{\end{eqnarray*}}
\newcommand{\bce}{\begin{center}}
\newcommand{\ece}{\end{center}}
\newcommand{\bpi}{\begin{picture}}
\newcommand{\epi}{\end{picture}}
\newcommand{\bfi}{\begin{figure} \begin{center}}
\newcommand{\efi}{\end{center} \end{figure}}

\newcommand{\bsl}{\begin{slide}{}}
\newcommand{\esl}{\end{slide}}

\newcommand{\bib}{thebibliography}

\newcommand{\ul}{\underline}

\newcommand{\sbs}{\subset}

\newcommand{\al}{\alpha}
\newcommand{\be}{\beta}
\newcommand{\ga}{\gamma}
\newcommand{\de}{\delta}
\newcommand{\ep}{\epsilon}

\newcommand{\la}{\lambda}

\newcommand{\si}{\sigma}

\newcommand{\De}{\Delta}

\newcommand{\bZ}{{\bf Z}}
\newcommand{\C}{{\Bbb C}}
\newcommand{\N}{{\Bbb N}}
\newcommand{\PP}{{\Bbb P}}
\newcommand{\R}{{\Bbb R}}
\newcommand{\Z}{{\Bbb Z}}

\newcommand{\cB}{{\cal B}}

\def\Zoverseven{\mathop{\lower 10pt \hbox{$ 
{\buildrel{\displaystyle\bar{Z}^2_{49}}\over{\hspace{-.2cm}
\scriptstyle{(7)}}} $}}  \nolimits}

\def\zovera{\mathop{\lower 10pt \hbox{$ {\buildrel{\displaystyle\bar{z}}
 \over {\scriptstyle{(a)}}} $}}{\lower 4pt
 \hbox{${\scriptstyle{ij}}$}}} 

\def\zoverk{\mathop{\lower 10pt \hbox{$ {\buildrel{\displaystyle\bar{Z}}
 \over{\scriptstyle{(k)}}} $}}
{\lower 4pt \hbox{${\scriptstyle{ij}}$}}}

\def\lczoverk{\mathop{\lower 10pt \hbox{$ {\buildrel{\displaystyle\bar{z}}
 \over{\scriptstyle{(k)}}} $}}
{\lower 4pt \hbox{${\scriptstyle{ij}}$}}} 

\def\Zovera{\mathop{\lower 10pt \hbox{$ {\buildrel{\displaystyle\bar{Z}} 
\over {\scriptstyle{(a)}}} $}}{\lower 4pt \hbox{${\scriptstyle{ij}}$}}} 
\def\zoverka{\mathop{\lower 10pt \hbox{$ {\buildrel{\displaystyle\bar{z}}
 \over
        {\scriptstyle{(k)}}} $}}{\lower 4pt \hbox{${\scriptstyle{m \ k+1}}$}}} 
\def\zoverkb{\mathop{\lower 10pt \hbox{$ {\buildrel{\displaystyle\bar{z} }
 \over
        {\scriptstyle{(k)}}} $}}{\lower 4pt \hbox{${\scriptstyle{k+1 \ m}}$}}} 
\def\zoverkc{\mathop{\lower 10pt \hbox{$ {\buildrel{\displaystyle\bar{z} } \over
        {\scriptstyle{(k)}}} $}}{\lower 4pt \hbox{${\scriptstyle{n \ k+2}}$}}} 
\def\zoverkd{\mathop{\lower 10pt \hbox{$ {\buildrel{\displaystyle\bar{z} } \over
        {\scriptstyle{(k)}}} $}}{\lower 4pt
 \hbox{${\scriptstyle{k-1 \ k+3}}$}}} 
\def\zoverke{\mathop{\lower 10pt \hbox{$ {\buildrel{\displaystyle z } \over
        {\scriptstyle{(k)}}} $}}{\lower 4pt 
\hbox{${\scriptstyle{k-1 \ k+2}}$}}} 
\def\zoverkf{\mathop{\lower 10pt \hbox{$ {\buildrel{\displaystyle\bar{ z}} \over
        {\scriptstyle{(k)(k+2)}}} $}}{\lower 4pt \hbox{${\scriptstyle{k-1 \ k+4}}$}}} 
\def\zoverkg{\mathop{\lower 10pt \hbox{$ {\buildrel{\displaystyle \bar{z}} \over
        {\scriptstyle{(k+2)}}} $}}{\lower 4pt \hbox{${\scriptstyle{k-2 \ k+1}}$}}} 
\def\zoverkh{\mathop{\lower 10pt \hbox{$ {\buildrel{\displaystyle \bar{z}} \over
        {\scriptstyle{(k)}}} $}}{\lower 4pt \hbox{${\scriptstyle{k-2 \ n-3+{i \over 2}}}$}}} 
\def\zoverki{\mathop{\lower 10pt \hbox{$ {\buildrel{\displaystyle z} \over
        {\scriptstyle{(k)}}} $}}{\lower 4pt \hbox{${\scriptstyle{k-1 \ k+2}}$}}} 
\def\zoverkj{\mathop{\lower 10pt \hbox{$ {\buildrel{\displaystyle z} \over
        {\scriptstyle{(k)}}} $}}{\lower 4pt \hbox{${\scriptstyle{k-1 \ k+2}}$}}} 
\def\vp{\varphi}
\def\g{\gamma}
\def \G{\Gamma}

\begin{document}
\title{\textsc{Braid Monodromy of Special Curves}
\author{Meirav Amram and Mina Teicher}}
 
\date{\today \\[1in]}

\maketitle
  
\begin{abstract}
In this article, we compute the braid monodromy of two algebraic curves
defined over $\R$. These two curves are of complex level not bigger than $6$,
 and they are unions of lines and conics. 
We use two different techniques for computing their braid monodromies.
These results will be applied to computations of fundamental groups of 
their complements in $\C ^2$ and $\C \PP^2$.  
\end{abstract}

\newpage

\noindent
{\Large {\bf \underline{Preface:}}}

\medskip

In this article we compute the braid monodromy of two curves: $S_1$ is a curve
 with two conic sections and $S_2$ is a curve with three conic sections.

We compute the braid monodromy by local computations in a neighbourhood of a 
single singular point, and by a series of Lefschetz diffeomorphisms on 
semicircles (see Theorem 2.25 and Theorem \ref{mon}).

The article is divided into 4 parts:
\begin{enumerate}
\item Introduction:  definition of braid monodromy.
\item An algorithm for computing the braid monodromy of line arrangements 
  and curves of degree $2$ and of complex level not bigger than $6$:
   \begin{itemize}
   \item[(2.1)] Definitions for $(\tilde D,K)$ - a model for $(D,K(x))$ 
        for all $x$ on the $x$-axis and definitions for $\de _{x_j}$ - 
        diffeomorphisms from one model to another.
   \item[(2.2)] $\xi _{x'_j} $ - skeletons in the fibres of 
        each $x$ and the corresponding skeletons $\la _{x_j}$ in the models.
   \item[(2.3)] The family of diffeomorphisms $\{ \be _x \}$ 
        from $(D,K(x))$ to $(\tilde D,K)$.
   \item[(2.4)] The Lefschetz diffeomorphisms $\Psi _h$ for 
        a semicircle $h$ below the real line with a radius $\al$.
   \item[(2.5)] Algorithm for computing the braid monodromy.
   \item[(2.6)] Notations of paths and the corresponding 
        halftwists.
   \item[(2.7)] Some examples for computing the action of $\de _{x_j}$ 
        on skeletons in $(\tilde  D,K)$.   
   \item[(2.8)] Motivation for having the two curves $S_1,S_2$ - 
        degeneration and regeneration.
   \end{itemize}
\item The braid monodromy of a curve $S_1$ with two conic sections.
\item The braid monodromy of a curve $S_2$ with three conic sections. 
\end{enumerate}
\subsection*{{\bf \underline{Notation:}}}
\begin{list}{}{\listparindent 0in \itemindent 0in \leftmargin \rightmargin}
\item $S_1, S_2$ algebraic curves defined over $\R\ ;\ S_1, S_2 \subset \C^2$.
\item $\pi_1 : S_i \rightarrow \  \C $ projection on the first 
coordinate (projection to x-axis) for $i = 1,2$.
\item $\pi_2 : S_i \rightarrow \  \C $ projection
 on the second coordinate (projection to y-axis) for $i = 1,2$.
\item$E(\mbox{ resp}.D)$ be a closed disk on the 
x-axis (resp. y-axis) with the center on the real
 part of the x-axis (resp. y-axis), s.t. \{singularities  of $ \pi_1\} \subseteq E \times (D-\partial D)$.
\item  $K(x) = \pi_2 (\pi^{-1}_1(x))$.
\item $K_\R(x) = K(x) \cap \R$.
\item $N = \{x \in E \mid \#K_\R(x) < n \}$.
\item $M$ a real number,  s.t. $x << M \;\;\;  \forall  x \in N$.
\item $E_\R = E \cap \R$.
\item $D_\R = D \cap \R$.
\item $\vp_M$ = the braid monodromy of an algebraic curve $S$ in $M$.
\item $B_n[D,K]$ = the braid group.
\item $\pi_1(\C^2 - S_i, u_0)$ = the fundamental group of a complement of a 
branch curve $S_i, \; i = 1,2$.
\item Paths will be signed by small letters and halftwists will be signed by capital letters.
\item We sign the points of $K(M)$ by $q_i \mbox{ or } q_j$.
\item $\underline{z}_{ij}$ = a path from $q_i \mbox{ to } q_j$ below the real line.
\item $\bar{z}_{ij}$ = a path from $q_i \mbox{ to } q_j$ above the real line.
\item $\stackrel{(a)}{\underline{z}}_{ij}$ = a path from $q_i \mbox{ to } q_j$, s.t. the path passes above $a$ and below the real line.
\item $\zovera$ = a path from $q_i \mbox{ to } q_j$, s.t. the path passes below $a$ and above the real line.
\item The corresponding halftwists: $H(\underline{z}_{ij}) = \underline{Z}_{ij}, H(\bar{z}_{ij}) = \bar{Z}_{ij}$.
\end{list}

\section{\underline{ Introduction: definition of braid monodromy}}

In this chapter, we define the braid monodromy and we present basic examples. 
To be able to define the braid monodromy, 
we have to enclose geometric model of the braid 
group and definitions of halftwist and g-base.

\bde {\bf Braid group $\mbox{\boldmath $B_n [D,K]$.}$}\\
Let $D$ be a topological disk, $K \sbs D$ be a finite set and
let $u \in\partial D$.
Let $\cB$ be the group of all diffeomorphisms $\be$ of $D$, s.t. $\be(K)=K$,
$\be|_{\partial D} = Id|_{\partial D}$. Such diffeomorphism acts naturally on
$\pi_1(D-K,u)$. We say that two such diffeomorphisms
are equivalent if they define the same automorphism on $\pi_1(D-K,u)$.
The quotient of $\cB$ by this equivalence relation is called {\bf the 
braid group $\mbox{\boldmath $B_n [D,K]$}$,} ($n=\#K$).

\ede
\bre {\rm 
We can define $B_n [D,K]$ as follows: there exists a natural 
$\psi :  \cB \to {\rm Aut} [\pi _1 (D-K, u)]$ and now we can define: 
$B_n[D,K]={\rm Im} (\psi)$, see [A]. Moreover, $B_n[D,K] \simeq B_n[D',K']$ for $(D,K)$ and $(D',K')$ topological disks, s.t.
$\# K=\# K'$.
(See [MoTe1], Lemma III.1.1).}
 \ere

\bde $\mbox{\boldmath $B_n$}$. \\
We can define $\mbox{\boldmath $B_n$}$ to be $B_n [D,K]$ for some $D,K$ as above. It is well defined 
by Remark 1.2.
\ede

\bde {\bf Halftwist of the closed interval $\mbox{\boldmath $[-1,1]$.}$} \\
Let $K=\{ -1,1 \} , \ \ \si=[-1,1]$. There exists a smooth and monotone function $\al$,
$\al: [0,2] \to [0,1]$, s.t.
$\al (x) =  \left\{ \matrix{
            1  &  x \in [0, {3 \over 2} ] \cr
            0  &  x \ge 2
} \right.
$.\\
Let us  look at $h: \{ x \ | \ |x| \leq 2 \} \to \{ x \ | \ |x| \leq 2 \} $, which is defined by:
$h(re^{i \varphi})=re^{i(\varphi +\al(r)\pi)},\  {\rm for}\  re^{i \varphi}=z \in \C .$
$h$ is called {\bf a halftwist w.r.t. $\mbox{\boldmath $[-1,1]$}$.}
\ede
\bre{\rm
In the domain $\{ z \in \C \ | \ |z| \leq {3 \over 2} \}$, $h(z)$ is a $180 ^{\circ} $ 
rotation counterclockwise. $h(z)={\rm Id}$ on $\{ z \in \C \ | \ |z|=2 \} $.
} \ere
\bde {\bf Halftwist w.r.t. a path from $a$ to $b$ ($H(\si)$).} \\
Let $D$ be a topological disk, $K \sbs D$, $\# K < \infty$, and let $a,b \in K$. Let 
$\si$ be a path from $a$ to $b$. Let $U$ be a closed neighbourhood  of $\si$, 
s.t. $U \cap K = \{ a,b \}$. Take $f: U \to \{ z \in \C \ | \ |z| \leq 2 \}$,s.t. 
$f(\si)=[-1,1]$. Look at the diffeomorphism defined by $(f \cdot h \cdot f^{-1} ) | _U$ 
and the identity map outside $U$. $\mbox{\boldmath $H(\si)$}$ is the braid defined by 
this diffeomorphism. 
\ede
\bde {\bf Motion from $K$ to $K'$.} \\
Let $K= \{ a_1, \cdots , a_n \}, \ K'= \{ a'_1, \cdots , a'_n \}$ be two sets and 
$K,K' \sbs {\rm Int} (D)$. {\bf A motion from $K$ to $K'$} inside $D$ is defined by $n$ 
continuous functions $m_i : I \to {\rm Int} (D), \ 1 \leq i \leq n$, which satisfy: 
$ m_i (0)=a_i$ and $m_i (1)=a'_i \ \ \forall i$; 
 $m_i (t) \not = m_j (t) \ \ \forall i,j \;\; , \;\;  i \not = j \;\; , \;\; 
\forall t \in I$.
\ede
\bre{\rm (cf. [A]).\\
If $K',K'' \sbs D$, then the motion from $K'$ to $K''$ inside $D$ induces a 
diffeomorphism up to isotopy of $D$ relative to the boundary from $(D,K')$ to $(D,K'')$. \\
A motion $m$ from $K$ to $K$ naturally induces a braid $\be _m \in B_n[D,K]$.
By this braid we can define the braid monodromy.
} \ere
\medskip

For defining the braid monodromy we look at the following situation: let $E$
(resp. $D$) be a closed disk on the $x$-axis (resp. $y$-axis) with the center 
on the real part of the $x$-axis (resp. $y$-axis).
Look at the real part of an 
algebraic curve $C$ included in $E \times D$.

Let $\pi _1: E \times D \to E, \ \pi_2: E \times D \to D$ be projections, and let 
$\pi =\pi _1 | _C : C \to E$, $\deg \pi =n$.   
Take $K(x)=\pi ^{-1} (x)$.

Let $N= \{ x \in E_\R \ | \ \# K_\R(x) \lvertneqq n \}$. Assume that $N \cap \partial E = \emptyset$. 
Choose $M \in \partial E$ and denote $K=K(M)=\pi ^{-1} (M)$, $K = \{ a_1, \cdots , a_n \}$.

In such a situation, we are going to introduce the notion of {\it braid monodromy}.
\bde {\bf Braid monodromy of a curve $\mbox{\boldmath$ C$}$ w.r.t. $\mbox{\boldmath $(E \times D, \pi _1,M)$.}$} \\
Let $M \in \partial E$, s.t. $Re(x) < Re(M), \ \forall x \in N$. Let $\si$ be a closed path 
in $E-N$. To each closed path in $E-N$ which begins and ends in $M$, 
there exist $n$ liftings to the curve. 
Each one of them begins and ends in points of $K(M)$, but the liftings 
are not necessarily closed. We project these liftings to the vertical fibre 
$\pi _1 ^{-1} (M)$ and we get $n$ paths in $\pi _1 ^{-1} (M)$, s.t. that each 
one of them begins and ends in points of $K(M)$. These paths induce 
a motion on $(D,K)$. This motion defines a braid in $B_n [M \times D , K]$. 
This braid is called  $\varphi ([\si])$. $\varphi$ is a homomorphism defined by:
$$\varphi  _M : \pi _1 (E-N,M) \to B_n [M \times D,K].$$
$\varphi _M$ is called {\bf the  braid monodromy of a curve $\mbox{\boldmath$ C$}$ w.r.t. 
$\mbox{\boldmath $(E \times D,\pi _1, M)$}$.  }    
\ede  

\bre{\rm See  examples in [MoTe2]. 
} \ere
\bde $\mbox{\boldmath $\ell (\ga)$.}$ \\
Let $u_0 \in \partial D$ and let $\ga$ be a simple path which connects $u_0$ to
a point $u \in K$, s.t. $\ga \cap (K-\{ u \} ) = \emptyset$. Let $U$ be 
a small closed disk centered at $u$ and let $c = \partial U$, 
s.t. $(K-\{ u \} ) \cap  {\rm Int} (U) =  \emptyset$. 
Then:
$\mbox{\boldmath $ \ell (\ga)$}$  $ = (\ga') ^{-1} c \ga ' \in \pi _1(D-K, u_0)$.
\ede
\input epsf
\begin{figure}[t]
\begin{center}
\setlength {\epsfysize}{1.8cm}
\epsfbox {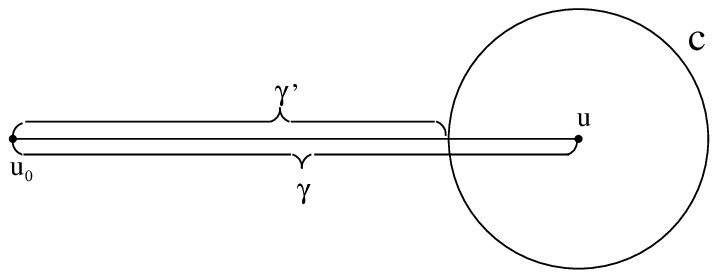}
\caption[]{$\ell (\ga)$}
\end{center}
\end{figure}

\bde {\bf A Bush, geometric base (g-base).}\\
Let $D$ be a disk, $K \subset D, \#K < \infty$.  Let $u_0 \in D-K$.
A set of simple paths $\{\g_i\}$ is  {\bf a Bush}  in $(D,K,u_0)$ if \ \ 
$\forall i, j \;\; \g_i \cap \g_j = u_0 ; \ \forall i \; \; \g_i \cap K =$ one point; $\g_i$ are ordered counterclockwise around $u_0$.
Let $\G_i = \ell(\g_i) \in \pi_1(D-K, u_0)$ be a loop around $K \cap \g_i$ determined by $\g_i$. $\{\G_i\}$ is called a {\bf a g-base} of 
 $\pi _1 (D-K,u_0)$.
\ede
\begin{figure}[htp]
\begin{center}
\epsfbox {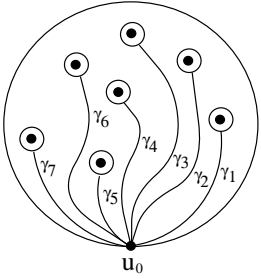}
\caption[]{A g-base}
\end{center}
\end{figure}

For computing the braid monodromy, we give the Emil Artin Theorem: 
\label{1.13}
\bthm {\bf Emil Artin Theorem}. \\
Let $C$ be a curve and let $\ell _1, \cdots, \ell _p$ be a g-base of $\pi _1 (\C -N, u)$.
Assume that the singularities of $C$ are: cusps, nodes, tangent points of a parabola/hyperbola with a line and branch points. \\
Then: for all $i$, there exist a half-twist $V_i$ and $r_i \in \Z$, s.t. 
$\varphi _M (\de _i)  = V_i ^{r_i}$
 and $r_i$ depends on the type of the singularity:
$$
r_i =  \left\{ \matrix{
            1  &  branch\ point  \cr
            2  &  node          \cr
            3  &  cusp          \cr
            4  &  tangent\ point         
} \right.
$$
\ethm

\vspace{-.3cm}

\noindent
{\em Proof:} See [A].
\hfill
$\Box$

We take curves which have the following properties:\\
1. \ $\forall x \in N, \ \# [ \pi _1 ^{-1}(x) \cap ({\rm singular\ points\ of}\ \pi _1)] =1$.\\
2. \  Each singular point of $\pi _1$ is one of the following types:\\
$a_1$= a branch point, which can be described locally by the equation $y ^2 -x=0;$   
   $a_2$= a branch point, which can be described locally by the equation
                $y ^2 +x=0;$    
  $b$= a tangent point, which can be described locally by the equation
                $y(y - x^2)=0;$      
   $c$= an intersection point of some lines or conics (parabolas or
              hyperbolas).

\bre{\rm
We are going to compute the braid monodromy of a part of a curve which 
is included in a finite interval. We choose two curves: $S_1$ with two 
conic sections and $S_2$ with three conic sections.  
} \ere

\section{\underline{An algorithm for computing the braid monodromy of line} \\
\underline{arrangements and curves of degree  $2$  and of complex level}\\
\underline{not bigger than $6$}}

In this chapter, we present the algorithm for computing braid monodromy
of a curve of complex level not bigger than $6$. In sections 2.1-2.4, we
introduce certain objects, models and notations that enable us in section
2.5 to present the algorithm. In section 2.6, we present notations of 
different halftwists. In section 2.7, we present a dictionary of Lefschetz
diffeomorphisms (using section 2.6). In section 2.8, we enclose notations for 
our algorithm. In Chapters 3 and 4, we shall demonstrate the algorithm 
from section 2.5 for two complex curves.  
\bde {\bf Complex level of a curve}. \\
A curve $C$ of degree $2$ is called a curve of  {\bf a complex level} $2$ if the 
curve includes exactly one curve of degree 2 (conic section). 
A curve $C$ of degree $2$ is called a curve of {\bf a complex level} $4$ if the 
curve includes exactly two conic sections. 
A curve $C$ of degree $2$ is called a curve of {\bf  a complex level} $6$ if the 
curve includes exactly three conic sections. 
\ede
\bre{\rm
 In the works of Moishezon-Teicher, there are computations 
of the braid monodromy of a curve of degree $2$ with one conic section.
In [Zu], there is an example of a computation of the
braid monodromy of a curve with two hyperbolas. In this work we are going 
to show the results of computations of the braid monodromy of two curves: 
$S_1$ which has two parabolas and $S_2$ which has three hyperbolas.
} \ere
\subsection{\underline{Models for $K(x)$}}

For these computations, we have  to define models of $K(x)$, $x \not\in N$, 
as follows:
Let $\tilde D = \{ z \in \C \ | \ |z- {n+1 \over 2}| \leq {n+1 \over 2} \}$ be a disk. The model 
depends on the complex level, therefore we have four types of models for $K(x)$

\begin{figure}[htp]
\begin{center}
\epsfbox {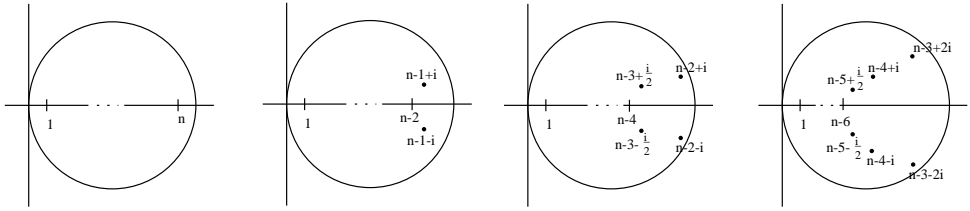}
\caption[]{Models for $K(x)$}
\end{center}
\end{figure}


\medskip

We have to define diffeomorphisms from each model to another too.
\bde $ \mbox{\boldmath $\De <k,l>: (\tilde D , K_i) \to (\tilde D, K_i), \ i=1,2,4,6$.}$\\ 
$\mbox{\boldmath $\De <k,l>$}$ is the diffeomorphism induced from the motion: $\{ k, \cdots , l \}$
are moving in a $180 ^{\circ}$ rotation counterclockwise centered at $k+l \over 2$.
\ede   

\label{2.4}
\bre{\rm
$\De <k,k+1>$ is actually the halftwist w.r.t. the 
line segment between $k$ and $k+1$.
} \ere
\label{2.5}
\bde $\mbox{\boldmath $\De ^{1 \over 2} _{I_2 \R} <k>: (\tilde D , K_2) \to (\tilde D, K_1)$.}$\\ 
$\mbox{\boldmath$\De ^{1 \over 2} _{I_2 \R} <k>$}$ $ = M_1 M_2 M_3, \ \ 1 \leq k \leq n-1$,
where:
$M_1=$ diffeomorphism induced from the motion: \ $n-1+i \to k + {1 \over 2} + {i \over 2}, \quad n-1-i \to k + {1 \over 2}- {i \over 2}$  ;  $M_2=$ diffeomorphism induced from the motion: 
$\{ k, \cdots ,n-2 \} \to \{ k+2 , \cdots , n \}$;
$M_3=$ diffeomorphism induced from the motion:  $ k + {1 \over 2} + {i \over 2} \to k , \quad k + {1 \over 2}- {i \over 2} \to k+1$
in a $90 ^{\circ}$ rotation counterclockwise.
\ede

\begin{figure}[htp]
\begin{center}
\epsfbox {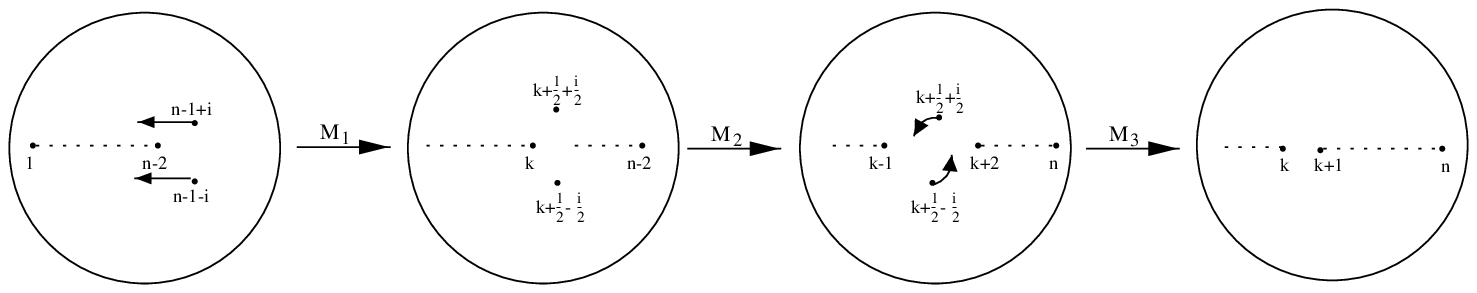}
\caption[]{$\De ^{1 \over 2} _{I_2 \R} <k>$}
\end{center}
\end{figure}


\bde $\mbox{\boldmath $\De ^{1 \over 2} _{\R I_2} <k>: (\tilde D , K_1) \to (\tilde D, K_2)$.}$ \\
$\mbox{\boldmath $\De ^{1 \over 2} _{ \R I_2} <k>$}$ $ = M_3 M_2^{-1} M_1^{-1}, \ \ 1 \leq k \leq n-1$,
where:
$M_3=$ diffeomorphism induced from the motion:
$k,k+1 \to k+ {1 \over 2} \pm {i \over 2}$
in a $90 ^{\circ}$ rotation counterclockwise;
$M_2^{-1}=$ diffeomorphism induced from the motion:
$\{ k+2, \cdots ,n \} \to \{ k , \cdots , n-2 \}$;
$M_1 ^{-1}=$ diffeomorphism induced from the motion:
 $ k + {1 \over 2} \pm {i \over 2} \to n-1 \pm i.$ 
\ede 
\begin{figure}[htp]
\begin{center}
\epsfbox {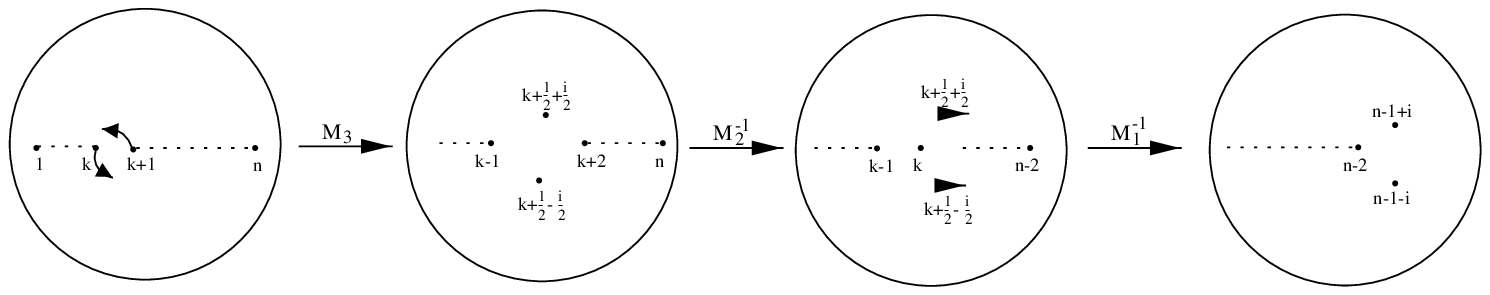}
\caption[]{$\De ^{1 \over 2} _{\R I_2} <k>$}
\end{center}
\end{figure}


\bde $\mbox{\boldmath $\De^{\frac{1}{2}'} _{I_4 I_2}  <k>: 
(\tilde D , K_4) \to (\tilde D, K_2)$.}$\\
$\mbox{\boldmath $\De^{\frac{1}{2}'} _{I_4 I_2} <k>$}$ $ = \hat{M_1} \hat{M_2} \hat{M_3}, \ \ 1 \leq k \leq n-3,$
where:
$\hat{ M_1}=$ diffeomorphism induced from the motion:
$n-2 \pm i \to n-2 \pm {i \over 2} ,\;  n-3 \pm {i \over 2 } \to n-3 \pm i$
composed with a diffeomorphism induced from the motion:
$n-2 \pm {i \over 2} \to k + {1 \over 2} \pm {i \over 2}$;
$\hat{ M_2}=$ diffeomorphism induced from the motion:
$n-3 \pm i \to n-1 \pm i, \; \{ k, \cdots ,n-4 \} \to \{ k+2 , \cdots , n-2 \}$ ;
$\hat{ M_3}=$ diffeomorphism induced from the motion:
$ k + {1 \over 2} \pm {i \over 2} \to k,k+1$
in a $90 ^{\circ}$ rotation counterclockwise.
\ede 
\begin{figure}[h]
\begin{center}
\epsfbox {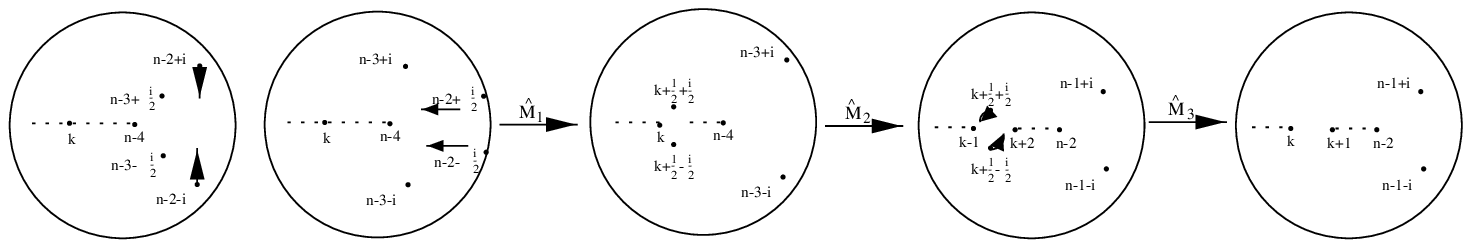}
\end{center}

\vspace{-4cm}

\caption{$\De^{\frac{1}{2}'} _{I_4 I_2} <k>$}
\end{figure}


\bde $\mbox{\boldmath $\De ^{1 \over 2} _{I_4 I_2} <k>: (\tilde D , K_4) 
\to (\tilde D, K_2)$.}$ \\ 
$\mbox{\boldmath $\De ^{1 \over 2} _{ I_4 I_2} <k>$}$ $ = \breve{M_1} \breve{M_2} \breve{M_3}, \ \ 1 \leq k \leq n-3,$
where:
$\breve{ M_1}=$ diffeomorphism induced from the motion:
$n-3 \pm {i \over 2} \to k + {1 \over 2} \pm {i \over 2}$;
$\breve{M_2}=$ diffeomorphism induced from the motion:
$n-2 \pm i \to n-1 \pm i, \; \{ k, \cdots ,n-4 \} \to \{ k+2 , 
\cdots , n-2 \}$;
$\breve{M_3}=$ diffeomorphism induced from the motion:
$ k + {1 \over 2} \pm {i \over 2} \to k,k+1$
in a $90 ^{\circ}$ rotation counterclockwise.
\ede 
\begin{figure}[htp]
\begin{center}
\epsfbox {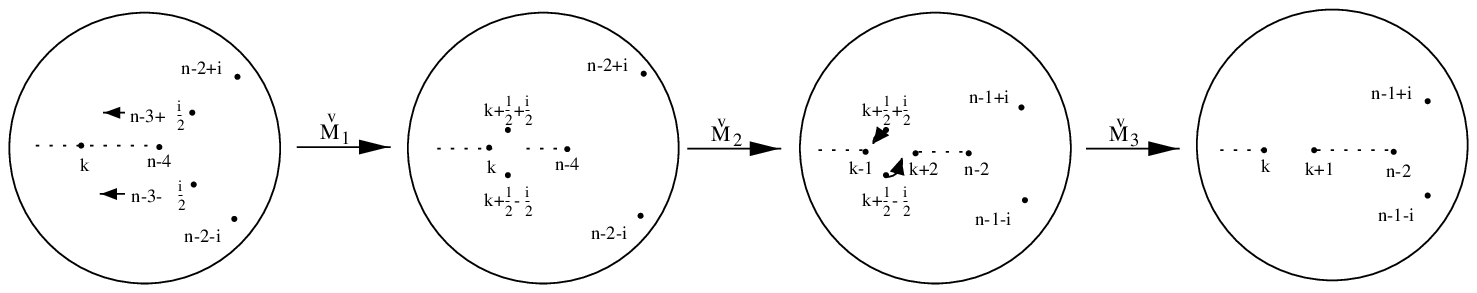}
\caption[]{$\De ^{1 \over 2} _{I_4 I_2} <k>$}
\end{center}
\end{figure}


\bde $\mbox{\boldmath$\De ^{1 \over 2} _{I_2 I_4} <k>: (\tilde D , K_2) 
\to (\tilde D, K_4)$.}$ \\
$\mbox{\boldmath$\De ^{1 \over 2} _{ I_2 I_4} <k>$}$ $ = \breve{M_3}
 \breve{M_2} ^{-1} \breve{M_1}^{-1}, \ \ 1 \leq k \leq n-3,$
where:
$\breve{ M_3}=$ diffeomorphism induced from the motion:
$k,k+1 \to k + {1 \over 2} \pm {i \over 2}$
in a $90 ^{\circ}$ rotation counterclockwise;
$\breve{M_2} ^{-1}=$ diffeomorphism induced from the motion:
$n-1 \pm i \to n-2 \pm i, \; \{ k+2, \cdots ,n-2 \} \to \{ k , \cdots ,
 n-4 \}$;
$\breve{M_1} ^{-1}=$ diffeomorphism induced from the motion:
$ k + {1 \over 2} \pm {i \over 2} \to n-3 \pm {i \over 2}.$
\ede

\begin{figure}[htp]
\begin{center}
\epsfbox {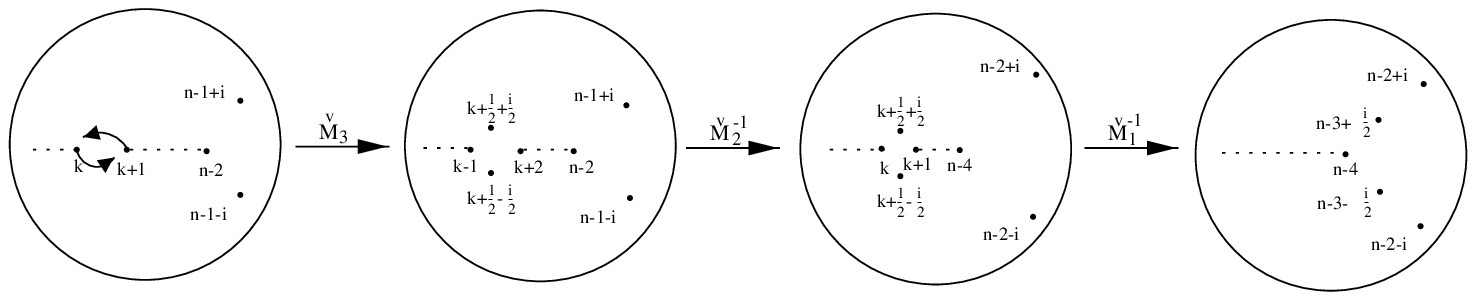}
\caption[]{$\De ^{1 \over 2} _{I_2 I_4} <k>$}
\end{center}

\vspace{-1cm}

\end{figure}


\bde $\mbox{\boldmath $\De  ^{1 \over 2} _{I_6 I_4} <k>:
 (\tilde D , K_6) \to (\tilde D, K_4)$. }$ \\
$\mbox{\boldmath $\De^{1 \over 2} _{ I_6 I_4} <k>$}$ 
$= \ddot{M_1} \ddot{M_2} \ddot{M_3}, \ \ 1 \leq k \leq n-5,$
where:
$\ddot{ M_1}=$ diffeomorphism induced from the motion:
$n-5 \pm {i \over 2} \to n-5 \pm i ,\;  n-3 \pm  2i  \to n-3 \pm 
{i \over 2}$
composed with a diffeomorphism induced from the motion:
$n-3 \pm {i \over 2} \to k + {1 \over 2} \pm {i \over 2}$;
$\ddot{ M_2}=$ diffeomorphism induced from the motion:
$n-4 \pm i \to n-2 \pm i, \; \{ k, \cdots ,n-6 \} \to \{ k+2 , \cdots , 
n-4 \}, \; n-5 \pm i \to n-3 \pm {i \over 2}$;
$\ddot{ M_3}=$ diffeomorphism induced from the motion:
$ k + {1 \over 2} \pm {i \over 2} \to k,k+1$
in a $90 ^{\circ}$ rotation counterclockwise.
\ede
\begin{figure}[htp]
\begin{center}
\epsfbox {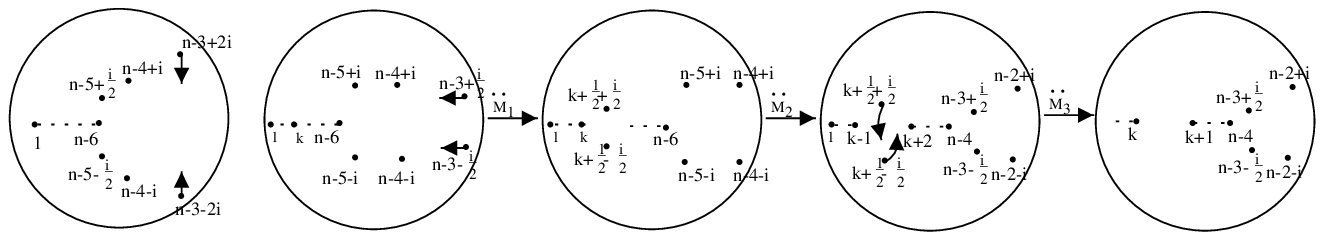}
\caption[]{$\De  ^{1 \over 2} _{I_6 I_4} <k>$ }
\end{center}
\end{figure}


\bde $\mbox{\boldmath $\De  ^{1 \over 2} _{I_4 I_6} <k>: 
(\tilde D , K_4) \to (\tilde D, K_6)$. }$ \\
$\mbox{\boldmath $\De  ^{1 \over 2} _{ I_4 I_6} <k>$}$ $
 = \ddot{M_3} \ddot{M_2}^{-1} \ddot{M_1}^{-1}, \ \ 1 \leq k \leq n-5,$
where:
$\ddot{ M_3}=$ diffeomorphism induced from the motion:
$ k,k+1 \to k + {1 \over 2} \pm {i \over 2}$
in a $90 ^{\circ}$ rotation counterclockwise;
$\ddot{ M_2}^{-1}=$ diffeomorphism induced from the motion:
$\{ k+2, \cdots ,n-4 \} \to \{ k , \cdots , n-6 \}, \; n-3 \pm 
{i \over 2} \to n-5 \pm i, \; n-2 \pm i \to n-4 \pm i $; $\ddot{ M_1}^{-1}=$ diffeomorphism induced from the motion:
$k + {1 \over 2} \pm {i \over 2} \to n-3 \pm {i \over 2}$ 
composed with a diffeomorphism induced from the motion:
$n-3 \pm {i \over 2} \to n-3 \pm 2i ,\;  n-5 \pm  i  \to n-5 
\pm {i \over 2}.$
\ede   

\begin{figure}[htp]
\begin{center}
\epsfbox {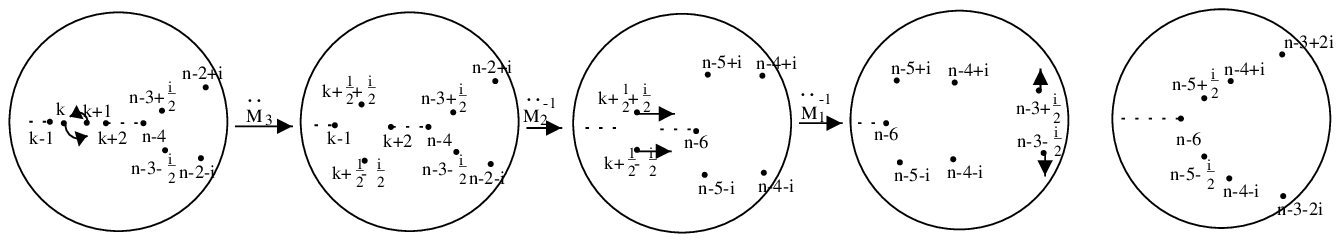}
\caption[]{$\De  ^{1 \over 2} _{I_4 I_6} <k>$}
\end{center}
\end{figure}


\vspace{-.3cm}

\subsection{\underline{Skeletons}}

Let $\pi = \pi_{1 \mid_C}, n = \deg \pi, N = \{x \in E \mid \# \pi^{-1} (x) 
< n \}, K_x = \pi^{-1}(x)$.  Therefore, for any real $x \notin N$, 
we have $n$ distinct
 real points $(x, y_i(x)), 1 \leq i \leq n$, in $K_x$.  
We choose a numeration
in $\{y_1(x) , \cdots , y_n(x)\}$, s.t. $y_1(x) < \cdots < y_n(x)$.

Let $\xi_x = \{\xi_{x,1} , \xi_{x,2} , \cdots , \xi_{x,n-1} \}$ be a 
sequence of real segments 
$[y_i(x), y_{i+1}(x)], 
1 \leq i \leq n-1$, in $x \times D$.

Now, we assume that $\forall x_j \in N$, there is only one singular
point of $C$ over $x_j$.

Let $x_j \in N$.  Choose $x'_j$, s.t. $\mid x'_j - x_j \mid = \alpha $,
$\alpha > 0$ a very small number.  Let $A_j$ be the singularity
of $C$ over $x_j$ and let $Y_j$ be the union of irreducible components
of $C$ containing $A_j$.  In $\{y_1(x'_j) , \cdots , y_n(x'_j) \}$, 
there is a sequence
with consecutive indices $\{y_{k_j}(x'_j) , y_{k_j+1}(x'_j) , \cdots ,
y_{\ell_j}(x'_j) \}$ which equals to $K_{x'_j} = Y_j \cap (x'_j \times D)$.

In this situation, we can define the following notation.

\bde $\mbox{\boldmath $(k_j, \ell_j)$ }$ {\bf Lefschetz pair of}
$\mbox{\boldmath $A_j$.}$ \\
The smallest and biggest indices $(k_{x_j}, \ell_{x_j})$ in the sequence 
considered above form a pair $(k_{x_j}, \ell_{x_j})$, which is called 
{\bf the Lefschetz pair of} $\mbox{\boldmath $A_j$.}$
\ede

Let $x_j \in N$ and let $A_j=(x_j,y_j)$ be a singular point above $x_j$, 
$A_j=\pi _1 ^{-1} (x_j)$.
Let $E_1 = \{ x \in E _{\R} -N \ | \ K(x)=n \}$. Take $x_j ' \in E_1$ 
and consider the following.

\bre  \ere
\begin{enumerate}
\item  If the point $A_j$ is of type $a_1, b,c$, we take $x_j ' = x_j+ \al$.
\item  If the point $A_j$ is of type $a_2$, we take $x_j ' = x_j-\al$.
\item  If the point $A_j$ is of type $a_1,a_2,b$, we take $l_{x_j}  = 
k_{x_j}+1$.
\end{enumerate}

\bre{\rm
Take $x_j, A_j, E_1, K(x'_j)$ as above.  We can determine the Lefschetz pair
$(k_{x_j}, \ell_{x_j})$ in a particular way:  Take a disk $D(A_j, \epsilon)$
for a radius $\epsilon$
 and a center at $A_j$.  So $\{y_k(x), \cdots , y_\ell
(x)\} \subset K(x'_j) \cap D(A_j, \epsilon)$.  This is a specific way
in order to determine 
$(k_{x_j}, \ell_{x_j})$, since here $D(A_j, \epsilon)$
is the disk containing $A_j$ which intersects $K(x'_j)$.
} \ere
We present an example.  
In the following curve $S$ $(\deg S = 4)$ there are three
singular points $A_1, A_2, A_3$.  Take their projections $x_1, x_2, x_3$ 
on the $x$-axis.
As above, consider $x'_1, x'_2, x'_3$.  Take $x'_j \times  
D$,  $j = 1,2,3$, the vertical 
line going up from $x'_j$ to the curve.  Draw
a disk $D(A_j, \epsilon)$ centered at $A_j$ for a very small radius
$\epsilon$.    We numerate all intersection points of $x'_j \times  D$
 with the
curve by an increasing order from 1 to 4. The Lefschetz pair is exactly the two
intersection and extreme points of 
$D(A_j, \epsilon)$ with the curve above $x'_j$.
\begin{figure}[h] 
\begin{center}
\setlength{\epsfysize}{4cm}
\epsfbox {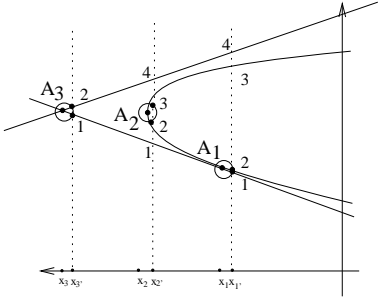}
\caption[]{Lefschetz pair}
\end{center}
\end{figure}


So, L-pair$(x_1) = (1,2),$ L-pair$(x_2)=(2,3)$, L-pair$(x_3)=(1,2)$.

We can not see more than 2 intersection points of $x'_3 \times D$ with the
curve since they are the intersection points of the conic with $x'_3 
\times D$, but they are complex.

In the same way as in the example we determine
 the L-pairs in the tables in Chapters 3,4.  In these tables 
we also consider a local numeration by going up on $x'_j \times D$.
For this reason we are going to
define skeletons and their models.

\bde {\bf Skeleton in  $\mbox{\boldmath $(D,K,K')$.}$}\\
Let $K' \sbs K \sbs D$ where $D$ is a disk and $K$ is a finite set.
 {\bf A skeleton
$\xi$ in $\mbox{\boldmath$(D,K,K')$}$} is a consecutive sequence of paths, each starting and ending 
in one of the points of $K'$. Each path starts where the previous one ends.
\ede
\vspace{-0.5cm}
\label{2.16}
\bre{\rm
We construct specific skeletons in $K(x_j ')$ \ $\forall x_j \in N$,
by considering the type of the singular point and the complex level.
Considering Figure 11, we determine the skeleton by
connecting the L-pair by a path.
} \ere

\bde  $\mbox{\boldmath $\xi  _{x_j '}$. }$ \\
Let $x_j \in N$, $A_j = (x_j, y_j)$ be a singular point of $\pi _1$. Let $\al, \ep$
be small enough in order to determine L-pair$(x_j)$. Let $x_j ' = x_j + 
\al$, 
L-pair $(x_j)$ =$(k,l)$. 
Let $\mbox{\boldmath$\xi  _{x_j '}$}$ be the following skeleton in $(D,K(x_j ') , K(x_j ') 
 \cap D(A_j,\ep))$ for $D(A_j, \epsilon)$ a disk around $A_j$ with a radius
$\epsilon$. 

\begin{enumerate}
\item For a point $A_j$ of type $a_1,b,c$: 
\begin{figure}[htp]
\begin{center}
\epsfbox {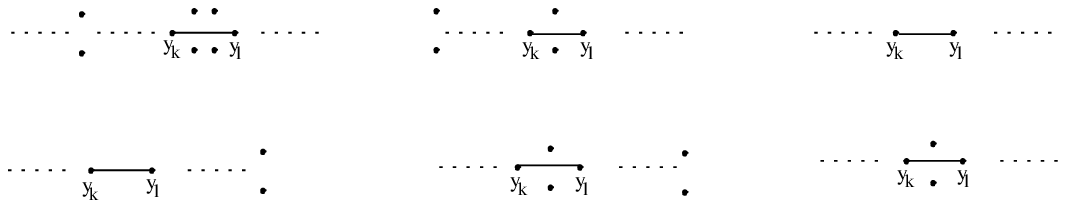}
\caption[]{$\xi  _{x_j '}$}
\end{center}
\end{figure}


\item For a point $A_j$ of type $a_2$:

\begin{figure}[htp]
\begin{center}
\epsfbox {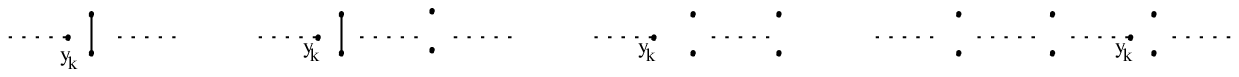}
\caption[]{$\xi  _{x_j '}$}
\end{center}
\end{figure}


\end{enumerate} 

\ede  
\vspace{-0.5cm}
\bre {\rm
In Figure 12 we present $\xi  _{x_j '}$ as a skeleton, connecting
the L-pair.  There are 4 possibilities for all points in all $x'_j \times 
D$.  It depends whether if in a specific $x'_j \times D$ the 
intersection points with conics are real or complex.  In Chapters 3,4, all
points can be real, 2 of them can be complex, 4 of them can be complex or
6 of them can be complex.

In Figure 13 we present skeletons which connect the two intersection 
points of a conic with a branch point of type $a_2$, with $x'_j \times 
D$.  Moving by Lefschetz diffeomorphism (see section 2.4), the
two points become complex.
} \ere

\subsection{\underline{Transition functions from $(D,K)$ to its models}}
For using the models, we have to define transition functions from $(D,K)$ to
$(\tilde D,K_i)$, $i=1,2,4,6$, and vice versa.

First, we define:
$E_1 = \{ x \in E_{\R} -N \ | \ \# K_{\R}(x)=n \}  \ ; \
E_2 = \{ x \in E_{\R} -N \ | \ \# K_{\R}(x)=n-2 \} \ ; \ 
E_4 = \{ x \in E_{\R} -N \ | \ \# K_{\R}(x)=n-4 \}\ ; \  
E_6 = \{ x \in E_{\R} -N \ | \ \# K_{\R}(x)=n-6 \} \ . \ $

These four sets satisfy:
$E_1 \cup E_2 \cup E_4 \cup E_6 = E_{\R}-N$.

\bre \ere 
\begin{enumerate} 
\item  $\forall x \in E_2$, there exist $z_1 (x),z_2(x) \in K(x)-K _{\R} (x)$, s.t. 
 ${\rm Re} (z_1 (x)) = {\rm Re} (z_2 (x)) $.
\item $\forall x \in E_4$, there exist $z_1 (x),z_2(x), z_3 (x),z_4(x) \in K(x)-K _{\R} (x)$, 
s.t.  ${\rm Re} (z_1 (x)) = {\rm Re} (z_2 (x)) $; ${\rm Re} (z_3 (x)) = {\rm Re} (z_4 (x))$.
\item $\forall x \in E_6$, there exist $z_i(x) \in K(x)-K _{\R} (x), \ 1 \leq i \leq 6$, 
s.t.  ${\rm Re} (z_i (x)) = {\rm Re} (z_{i+1} (x)), \ i=1,3,5$.
\end{enumerate} 

 We prefer to work in models, in which points are arranged in a more
comfortable way.  Now we define a family of diffeomorphisms 
from our disk $D$ to
$\tilde{D}$ its model.  Later we define also the models of 
the skeletons defined above. 

\ble
Let $L$ be a connected component of $E_{\R} -N$. Then, there exist $i$, 
$i=1,2,4,6$, s.t. $L \sbs E_i$. \\
Moreover, there exists a continuous family of diffeomorphisms 
$\{ \be _x \ | \ x \in L \}$, s.t.
\begin{enumerate}
\item $\be _x : D \to \tilde D$.  
\item $\be _x (K(x)) = K_i, \ \ i=1,2,4,6$. 
\item $\be _x (D_{\R}) = \tilde D \cap \R$. 
\item $\forall x,x' \in L, \ \forall y \in \partial D: \ \be _x (y)= \be _{x'} (y)$.  
\end{enumerate}
\ele

\noindent
{\it Proof:} Standard proof in complex analysis.
\hfill
$\Box$
\bre{\rm
$\be _x ^{\vee}$ is the isomorphism induced naturally from $\be _x$,
i.e.: if $\be _x$ satisfies the conditions of the lemma, then:
$\be _x ^{\vee}: B_n [D,K(x)] \to B_n [\tilde D, K_i]$
is an isomorphism, $ i=1,2,4,6$.
} \ere

We need to construct models for the skeletons. Let us look at the following
definition:
\bde  $\mbox{\boldmath $\la _{x_j}$}$ {\bf for} 
$\mbox{\boldmath $x_j \in N$.}$\\
Let $x_j \in N$, $A_j=(x_j,y_j)$ a singular point of $\pi _1$ above $x_j$. Let
L-pair$(x_j)=(k_{x_j},l_{x_j})$. Let $\mbox{\boldmath $\la _{x_j}$}$ 
be the following skeleton in $(\tilde{D}, K_i)$.
\ede
\begin{enumerate}
\item {\em for $A_j = a_1,b,c$.}
\begin{figure}[htp]
\begin{center}
\epsfbox {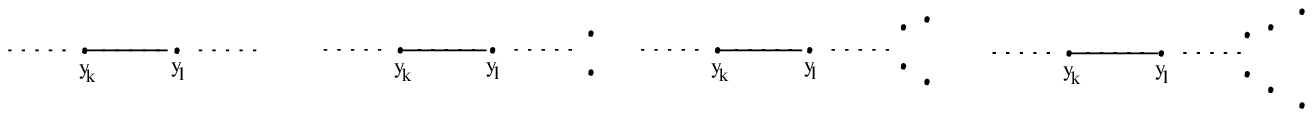}
\caption[]{$\la _{x_j}$}
\end{center}
\end{figure}


\item {\em for $A_j=a_2$.}
\begin{figure}[htp]
\begin{center}
\epsfbox {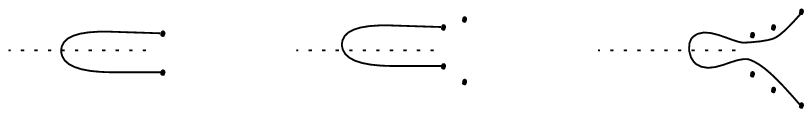}
\caption[]{$\la _{x_j}$}
\end{center}
\end{figure}

\end{enumerate}  

These skeletons are easier to work with, since complex points
are signed in the end (see Figures 14, 15).
 
\subsection{\underline{Lefschetz diffeomorphism induced by a curve $h(t)$}}
\label{diff}
Let $h(t)$, $t \in [0,1]$, be a curve in $E-N$. $h$ is lifted to $n$ paths in $C$ 
and projecting them to $D$, we get $n$ paths in $D$, each starting in one
of the points of $K(h(0))$ and ending in one of the points of $K(h(1))$. 
Thus, we get a motion $(D,K(h(0))) \to  (D,K(h(1)))$. This motion induces 
a homotopically equivalence class of diffeomorphisms, denoted $\Psi _h$.
\bde  $\mbox{\boldmath $\de _{x_j}$}$ {\bf for} 
$\mbox{\boldmath $x_j \in N$.}$\\
$\mbox{\boldmath $\de _{x_j} $}$ $= \De _{I_6 I_4} ^{1 
\over 2} < k_{x_j}>$ for $A_j=a_1$ and there 
is a transformation from a model with a complex level $6$ to a model with 
a complex level $4$, s.t. the two rightest complex points become real.
\\ 
$\mbox{\boldmath $\de _{x_j} $}$ $= 
\De _{I_4 I_2} ^{1 \over 2} < k_{x_j}>$ for $A_j=a_1$ and there 
is a transformation from a model with a complex level $4$ to a model with 
a complex level $2$, s.t. the two close complex points become
 real.\\ 
 $\mbox{\boldmath $\de _{x_j}$}$ $=  \De^{\frac{1}{2}'} _{I_4 I_2}
< k_{x_j}>$ for $A_j=a_1$ and there 
is a transformation from a model with a complex level $4$ to a model with 
a complex level $2$, s.t. the two righter  
complex points become real.\\ 
$\mbox{\boldmath $\de _{x_j} $} $ $ = \De _{I_2 \R} ^{1 \over 2} < k_{x_j}>$ for $A_j=a_1$ and there 
is a transformation from a model with a complex level $2$ to a real model,  
 s.t. the two complex points become real.\\ 
$\mbox{\boldmath $\de _{x_j} $} $ $ = \De _{\R I_2} ^{1 \over 2} < k_{x_j}>$ for $A_j=a_2$ and there 
is a transformation from a real model to a model with 
a complex level $2$, s.t. two real points become complex.\\
$\mbox{\boldmath $\de _{x_j}$} $ $ = \De _{I_2 I_4} ^{1 \over 2} < k_{x_j}>$ for $A_j=a_2$ and there 
is a transformation from a model with a complex level $2$ to a model with 
a complex level $4$, s.t. two real points become complex and the closer
ones.\\ 
$\mbox{\boldmath $\de _{x_j}$} $ $ = \De _{I_4 I_6} ^{1 \over 2} < k_{x_j}>$ for $A_j=a_2$ and there 
is a transformation from a model with a complex level $4$ to a model with 
a complex level $6$, s.t. two real points become complex and the 
rightest ones.\\ 
$\mbox{\boldmath $\de _{x_j} $} $ $ = \De ^2 <k_{x_j},l_{x_j}> \ $ for $A_j=b$.\\
$\mbox{\boldmath $\de _{x_j} $} $ $ = \De  <k_{x_j},l_{x_j}> \ $ for $A_j=c$.
\ede

\medskip

In the following theorem we give an algorithm for computing a Lefschetz 
diffeomorphism of a halftwist counterclockwise around a point in $N$ below the
real line.  
\bthm\label{beta_x0}
Let $x_j \in N$, where $N \subset E_\R$. Let $h$ be a semicircle of radius $\al$ from $x_j -\al$ to 
$x_j +\al$ below the real line. Let $\Psi _h$ be the Lefschetz diffeomorphism 
induced by $h$. Let L-pair$(x_j)= (k_{x_j},l_{x_j})$, $A_j = (x_j,y_j)$ a singular point of $\pi _1$. \\
Then: $\be ^{-1} _{x_j-\al} \cdot \Psi _h \cdot  \be _{x_j+\al}=\de _{x_j}$, where $ \de _{x_j}$ 
is defined above, considering the type of the singularity of $A_j$.
\ethm

\noindent
{\it Proof:} See [MoTe2], Proposition 1.2.
\hfill
$\Box$

In the following theorem we present an algorithm for computing 
the Lefschetz diffeomorphism of a curve below the real line. Remember that 
$N \sbs x$-axis.
\bthm 
Let $x_j \in N$. Let $\al,\ep$ be small enough, s.t. L-pair$(x_j)$ 
is well determined.
Let $x_j ' = x_j + \al$. Let $\ga _j '$ be a path from $x_j '$ to $M$ below 
the real line. \\
Then:
$(\xi _{x_j '}) \Psi _{\ga _j '} \sim (\la _{x_j})\displaystyle
 ({\prod _{m=j-1} ^1} \de _{x_m}) \be _M ^{-1} 
\qquad (\sim={\rm homotopic}).$
\ethm

\noindent
{\it Proof:} See [MoTe2], Proposition 1.3.
\hfill
$\Box$
\bde $\mbox{\boldmath ${\rm L.V.C.} (\ga_j)$}$.\\
Let $x_j \in N$. Let $\ga _j$ be a path from $x_j$ to $M$ below the real line, 
s.t. $x_j ' \in \ga _j$. Let $\ga _j '$ be the part of $\ga _j$ from $x_j '$ to $M$.
Let $\Psi _{\ga _j '}$ be the Lefschetz diffeomorphism corresponding to $\ga _j '$.\\
Define {\rm L.V.C.}  of $\ga _j$ by: 
$\mbox{\boldmath ${\rm L.V.C.} (\ga_j)$}$ = $(\xi _{x_j '}) \Psi _
{\ga _j '}$.
\ede

\subsection{\underline{Algorithm for computing the braid monodromy}}\label{com}
For all $x_j \in N$, we define $\ep _{x_j}$ w.r.t. $A_j =  (x_j,y_j)$ as follows:   
\bde $\mbox{\boldmath $\ep _{x_j}$.}$ \\
Let $A_j =  (x_j,y_j)$ be the singular point of $\pi _1$ above $x_j$. Define:
$$
\mbox{\boldmath $\ep _{x_j}$} =  \left\{ \matrix{
            1  &  A_j=a_1 \ {\rm or} \ A_j=a_2   \cr
            4  &  A_j=b          \cr
            2  &  A_j=c   \ .      
} \right.
$$
\ede
\bthm\label{mon}
Let $\varphi _M$ be the braid monodromy of a curve $C$ w.r.t. $(E \times D, \pi _1,M)$.
Let $x_j \in N$, $\forall j$, (where $N \subset E_\R$), and let $\ga _j$ be paths below the real line, 
which connect $M$ to $x_j$, $\forall j$. Then:
$\varphi _M (\ell (\ga _j)) = \De <(\xi _{x_j '})\Psi _{\ga _j '}>
 ^{\ep _{x_j}}\ , \ \ \forall j.$\\
Equivalently:
$\varphi _M (\ell (\ga _j)) = (Q ^{-1} \De (\la _{x_j}) Q)
 (\be _M ^{-1}) ^{\vee}\ \ , \ \  where \ \  \displaystyle Q= 
({\prod _{m=j-1} ^1} \de _{x_m})^{\vee}.$
\ethm

\noindent
{\it Proof:} See [MoTe2], Proposition 1.5.
\hfill
$\Box$
\bco\label{col-algo}
$\varphi _M (\ell (\ga _j)) = \De < {\rm L.V.C.}(\ga _j)> ^{\ep _{x_j}}= \De < (\la _{x_j})( \displaystyle {\prod _{m=j-1} ^1 \de _{x_m}})\be _M ^{-1}> ^{\ep _{x_j}}$.
\eco

\noindent
{\bf The algorithm:} \\
Let $x_j \in N$. Let $\ga _j$ be a path below the real line from $x_j$ to $M$. 
We want to compute $\varphi _M (\ell (\ga _j))$. Corollary 2.29 gives us 
an algorithm for computing the braid monodromy of $\ell (\ga _j)$ by 
L.V.C.$(\ga_j)$. L.V.C.$(\ga_j)$ is computed by the Lefschetz 
diffeomorphism of $\ga _j '$, which acts on $\xi _{x_j '}$ (see Theorem 2.25).

\medskip 

\bre{\rm
If we replace $\ga _j$ by a path above the real line, 
the formulas for $(\xi _{x_j '} ) \Psi _{\ga _j '}$ and for $\varphi _M (\ell (\ga _j))$ will be
replaced by formulas in which all diffeomorphisms are induced from clockwise
motions. If we replace $0 \ll M$ by $M \ll 0$, then we take clockwise motion 
in the different $\de _i$, we apply $\de _i$ in a reversed order, and also the
position of $a_1,a_2$ are exchanged. If we do both changes, then we have 
to change the order of the $\de _i$'s. Moreover, the positions of $a_1,a_2$ 
are exchanged.
} \ere
\subsection{\underline{Symbols for paths and halftwists}}\label{symb}
Paths will be signed by small letters  and halftwists will be signed 
by capital letters. We sign the points of $K(M)$ by $q_i$ or $q_j$.\\
$\ul{z} _{ij} = $ a path from $q_i$ to $q_j$ below the real line. \\  
$\bar{z} _{ij} = $ a path from $q_i$ to $q_j$ above the real line. \\  
$\stackrel{(a)}{\ul{z}} _{ij} =$ a path from $q_i$ to $q_j$, s.t. the path passes
above $a$ and below the real line.\\
$\zovera =$ a path from $q_i$ to $q_j$, s.t. the path passes
below $a$ and above the real line.

The paths can be written as conjugations, for example:\\
$\ul{z}^{\ul{Z}^2_{ik}}_{ij}$ for $i < k < j$ is: 
$\vcenter{\hbox{\epsfbox{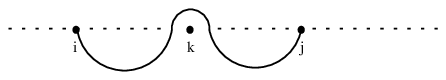}}}$, 
where $Z^2_{ik}$ is a halftwist
taken twice when $i$ moves around $k$ counterclockwise below the
real line.

Another example: 
$\vcenter{\hbox{\epsfbox{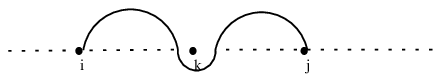}}}$
is $(\bar{z}_{ij})^{\bar{Z}^{-2}_{ik}}$, where
$i$ 
moves 
around $k$ clockwise above the real line in a 360$^0$ rotation.  
This 
path is also: $\lczoverk$.
Another example is:  $\vcenter{\hbox{\epsfbox{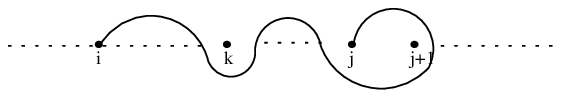}}}$.
This path is:$(\lczoverk)^{Z^2_{j \ j+1}}$, where
$j$ moves around $j+1$ counterclockwise in a 360$^0$ rotation.

In Chapters 3,4 we present the paths  which are more complicated  and may be
signed as a path conjugated by some halftwists, as explained above.

The corresponding halftwists are:
$H(\ul{z} _{ij}) = \ul{Z} _{ij} \ , \ 
H(\bar{z} _{ij})= \bar{Z} _{ij} \ , \ 
H(\stackrel{(a)}{\ul{z}} _{ij}) =\stackrel{(a)}{\ul{Z}} _{ij}\ , \  
H(\zovera) = \Zovera$.

\subsection{\underline{Examples for computing the
 action of $\de _{x_j}$ on skeletons
in 
$(\tilde D,K_i)$}}

Remember that $\De <k,k+1>$ is a $180 ^{\circ}$ rotation counterclockwise, 
while $k$ and $k+1$ are connected by a line: 
$\vcenter{\hbox{\epsfbox{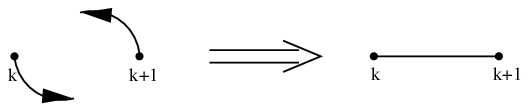}}}$

After the rotation, we continue to call the points by their places. If there
is $\De ^r <k,k+1>$ for $r \in \N$, then it is a  $180 ^{\circ}$ 
rotation counterclockwise, which is taken $r$ times. Hence:
\begin{enumerate}
\item If $r$ is even, then $k,k+1$ are coming back to their original places.  
\item If $r$ is odd, then $k,k+1$ are changing places. 
\end{enumerate} 

Now, we give some examples for all types of $\de _{x_j}$ which were defined in section \ref{diff}
 
\def\zoverkdt{\mathop{\lower 10pt \hbox{$ {\buildrel{\displaystyle\bar{z} } \over{\scriptstyle{(k)}}} $}}{\lower 4pt \hbox{${\scriptstyle{k-1 \ k+2}}$}}}
 
\def\zoverkx{\mathop{\lower 10pt \hbox{$ 
{\buildrel{\displaystyle\bar{z} } \over
        {\scriptstyle{(n-5 + \frac{i}{2})(n-4 + i)}}} $}}{\lower 4pt
 \hbox{${\scriptstyle{k \ n-3+2i}}$}}}

\bigskip
\noindent
$\De <k,k+1> (\bar{z} _{mk}) = \zoverka$.
$\hspace{.5cm}\vcenter{\hbox{\epsfbox{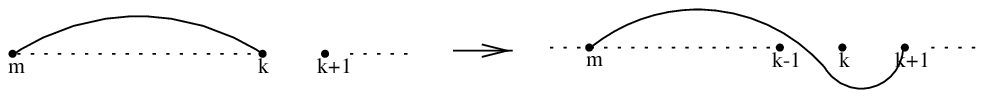}}}$\\

\bigskip

\noindent
$\De _{I_4 I_2} ^{1 \over 2} <k> (\ul{z} _{mn}) = 
\stackrel{(k+1)}{{\ul{z}}} _{m \ n+2}$.
$\hspace{.5cm}\vcenter{\hbox{\epsfbox{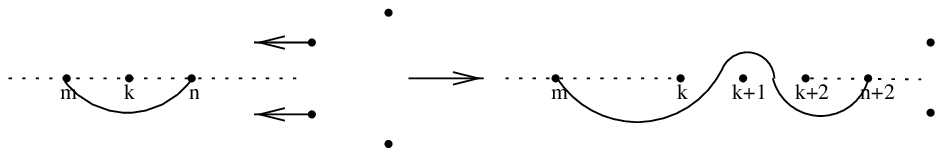}}}$

\bigskip

\noindent
$\De^{\frac{1}{2}'} _{I_4 I_2} <k>
 (\stackrel{(k+1)}{z} _{k-1 \ k+2})
= \zoverkf$.
$\hspace{.5cm}\vcenter{\hbox{\epsfbox{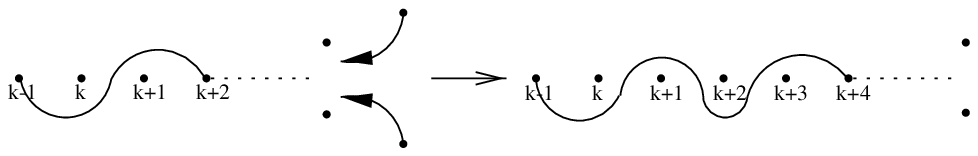}}}$

\bigskip

\noindent
${\De _{I_2 I_4} ^{1 \over 2}}  <k>  
(z _{k+2 \ k+3}^{Z^2_{k+1 \ k+2}})=
z_{k \ k+1}^{ \bar{Z}^{-2}_{k+1 \ n-3 + \frac{i}{2}} }$.
$\hspace{.5cm}\vcenter{\hbox{\epsfbox{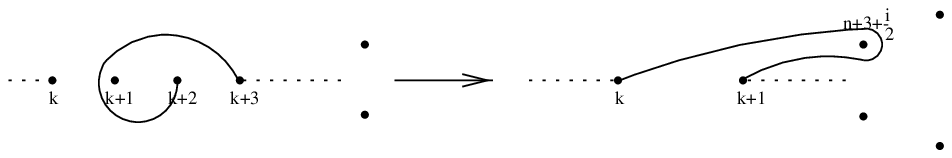}}}$

\bigskip
\bigskip

\noindent
${\De _{\R I_2} ^{1 \over 2}}  <k> (\stackrel{(k)}{z}_{k-1 \ k+2})= 
z _{k-1 \ k}$.
$\hspace{.5cm}\vcenter{\hbox{\epsfbox{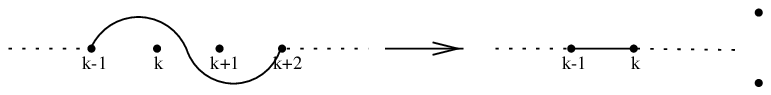}}}$

\bigskip
\bigskip
\noindent
${\De _{I_2 \R} ^{1 \over 2}}  <k>
 (z_{k-1 \ k}) = 
\stackrel{(k+1)}{\ul{z}} _{k-1 \ k+2}$.
$\hspace{.5cm}\vcenter{\hbox{\epsfbox{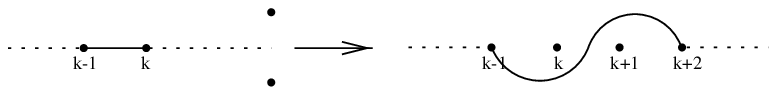}}}$

\bigskip
\bigskip
\noindent
$\De _{I_6 I_4} ^{\frac{1}{2}}  <k> (z_{k-1 \ k}) = 
\zoverkj$.
$\hspace{.5cm}\vcenter{\hbox{\epsfbox{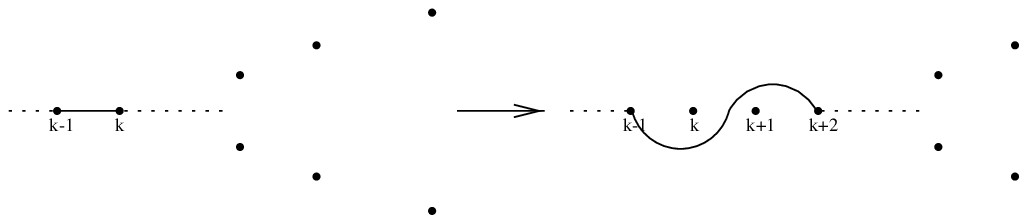}}}$

\bigskip

\noindent
$\De _{I_4 I_6} ^{\frac{1}{2}}  <k>  (z _{k+1 \ k+2})  
 = \zoverkx$.
$\hspace{.2cm}\vcenter{\hbox{\epsfbox{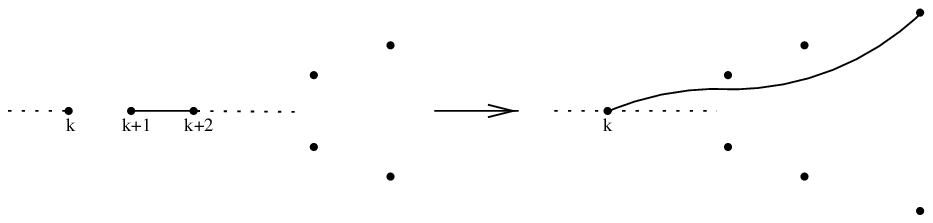}}}$
\bigskip

\label{2.31}
\bre{\rm $(\xi _{x_j '}) \Psi _{\ga _j '}$ is of type $z_{mk}$  ($m,k \in \C$), 
and $\varphi _M (\ell (\ga _j))$ is $Z_{mk} ^{\ep _{x_j}}$ the corresponding halftwist 
(with the corresponding power according to the singularity type.) 
} \ere

\subsection{\underline{Motivation for having the
 two curves $S_1,S_2$ - 
degeneration and} \\
\underline{regeneration}}

In this work we compute the braid monodromy of two curves $S_1,S_2$. These curves
appear in a regeneration process of a certain branch curve of an algebraic
surface.
\bde {\bf Degeneration}.\\
Let $X,Y$ be two projective varieties. We say that $X$ is {\bf a degeneration} 
of $Y$ if there exist an irreducible variety $V$ and an onto map
$\pi : V \to \C$ which satisfies: $\pi ^{-1} (0) \cong X, \pi ^{-1} (1) \cong Y$. 
$\pi ^{-1} (1)$ is a generic fibre. 
\ede

We are interested in a degeneration of a hypersurface to a union of planes.
Let us take $Y \sbs \C \PP ^3$ a surface of degree $5$, and we degenerate it 
to a union of planes $\hat Y \sbs \C \PP ^N$, s.t. each plane is $\C \PP ^2$. 
We take a generic projection of $Y$ to $\C \PP ^2$, and we get $C \sbs \C \PP ^2$
a branch curve. We take a generic projection of $\hat Y$ to $\C \PP ^2$ and 
we get a branch curve $\hat C \sbs \C \PP ^2$, which is a union of lines.

Regeneration is the opposite process of degeneration. So we regenerate $\hat C$
to get $C$:
\clearpage
\begin{figure}[tbh]
\begin{center}
\epsfbox {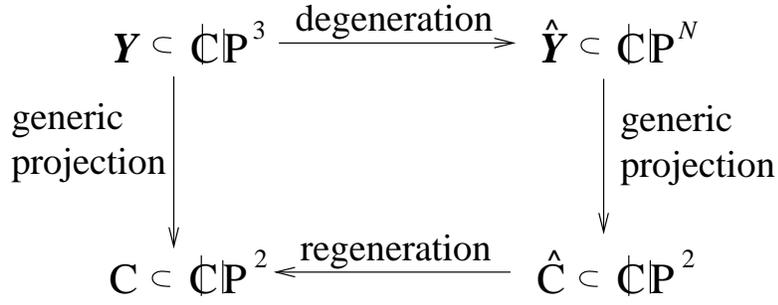}
\caption[]{A diagram}
\end{center}
\end{figure}


In their works, Moishezon-Teicher computed a regeneration of the degenerated
object in which six lines intersect. Here we regenerate 
an intersection point of six lines which do not appear in the works of 
Moishezon-Teicher. 

Regeneration process is opposite in nature to the degeneration process. The
braid monodromy of the regenerated curve depends very much on the type
and the order of the degeneration. In [MoTe4] Moishezon and Teicher computed
the braid monodromy of a certain 6-point regeneration. In this work, 
we compute braid monodromy of curves which will be needed in the future, 
and appear in the regeneration process of a different degeneration.  

During the process of the regeneration, we will get two possibilities:
\begin{enumerate}
\item The first regeneration produces a curve $S_1$, 
   which is treated in Chapter \ref{S1}
\item The second regeneration produces a curve $S_2$, 
   which is treated in Chapter \ref{S2}
\end{enumerate} 
The steps of the regeneration are:
\begin{figure}[htp]
\begin{center}
\epsfbox {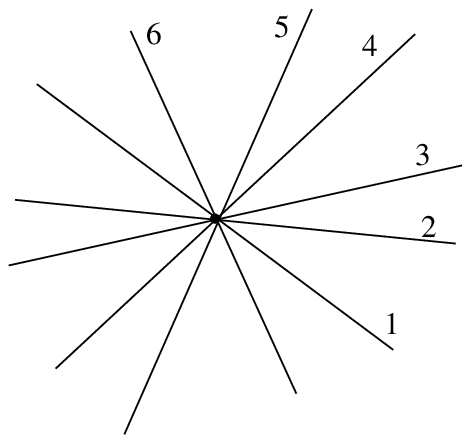}
\caption[]{Step 1}
\end{center}
\end{figure}
\clearpage

\begin{figure}[htp]
\begin{center}
\epsfbox {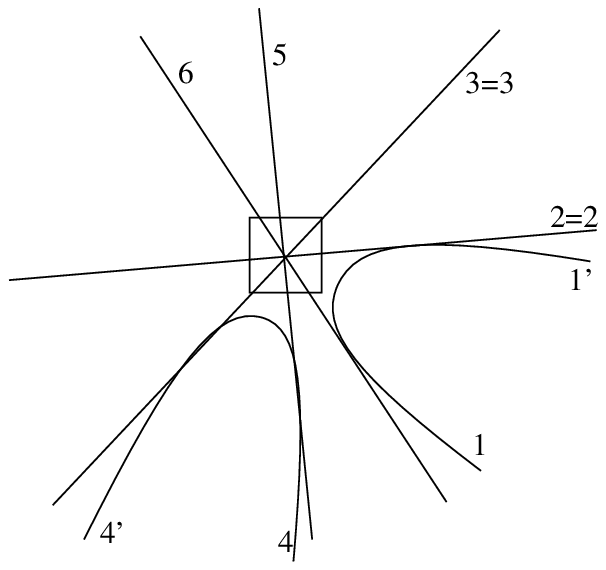}
\caption[]{Step 2}
\end{center}
\end{figure}


The regeneration of the lines $l_2,l_3,l_5,l_6$:
\begin{figure}[htp]
\begin{center}
\epsfbox {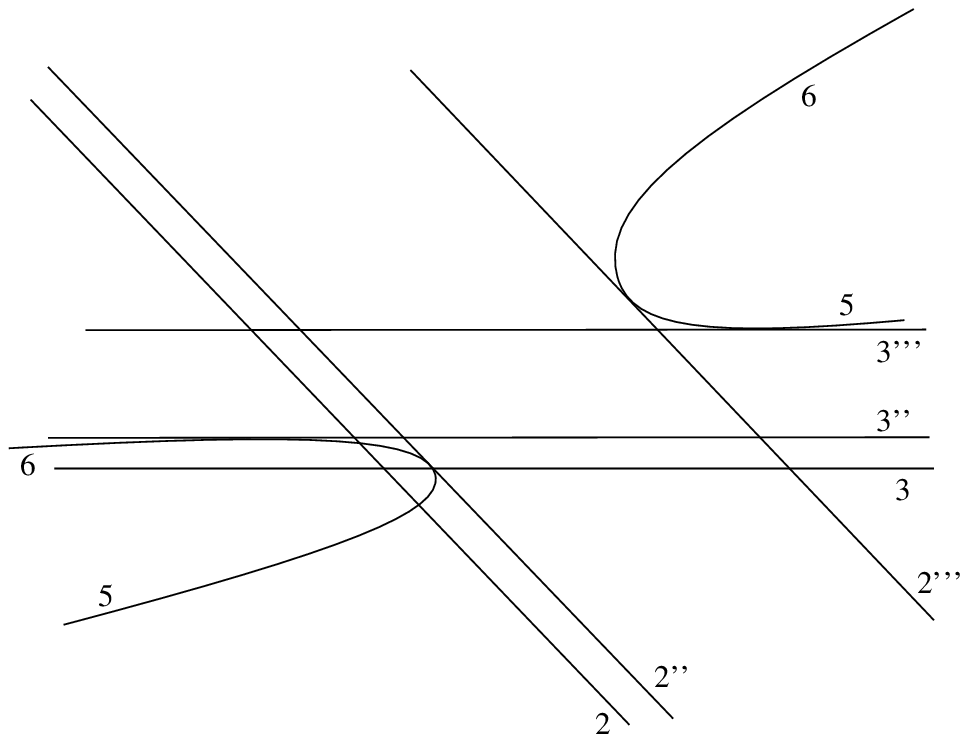}
\caption[]{Step 3}
\end{center}
\end{figure}

\clearpage
The two possibilities are $S_1$ in Figure 20  and $S_2$ in Figure 21:
\begin{figure}[htp]
\begin{center}
\setlength{\epsfysize}{6.8cm}
\epsfbox {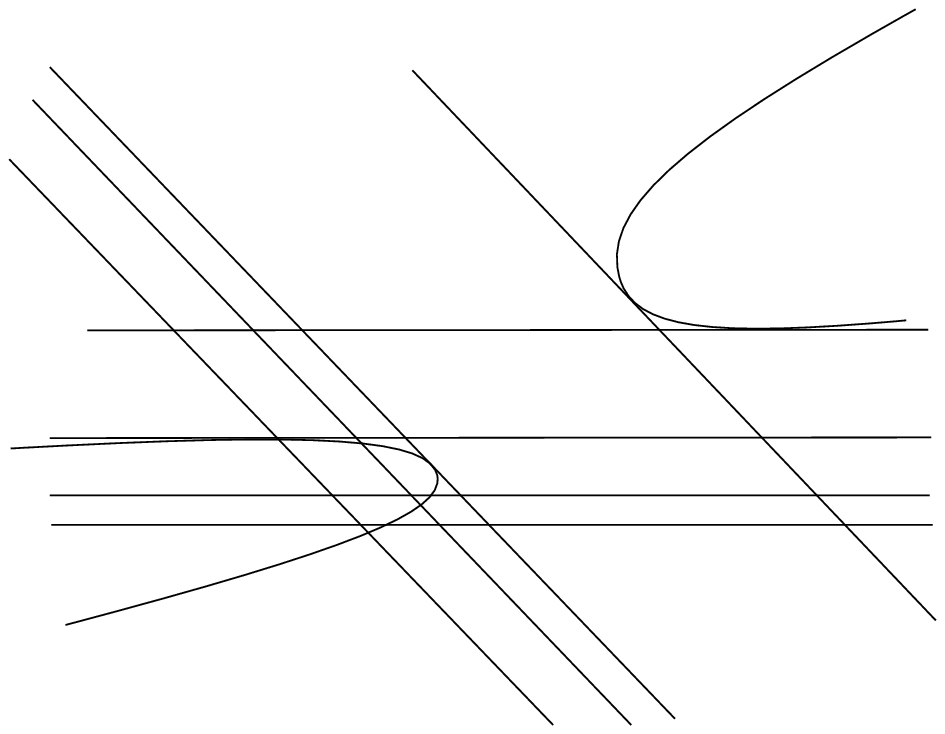}
\caption[]{$S_1$} \label{S2_fig}
\end{center}
\end{figure}
\vspace{2cm}

\begin{figure}[htp]
\begin{center}
\setlength{\epsfysize}{7cm}
\epsfbox {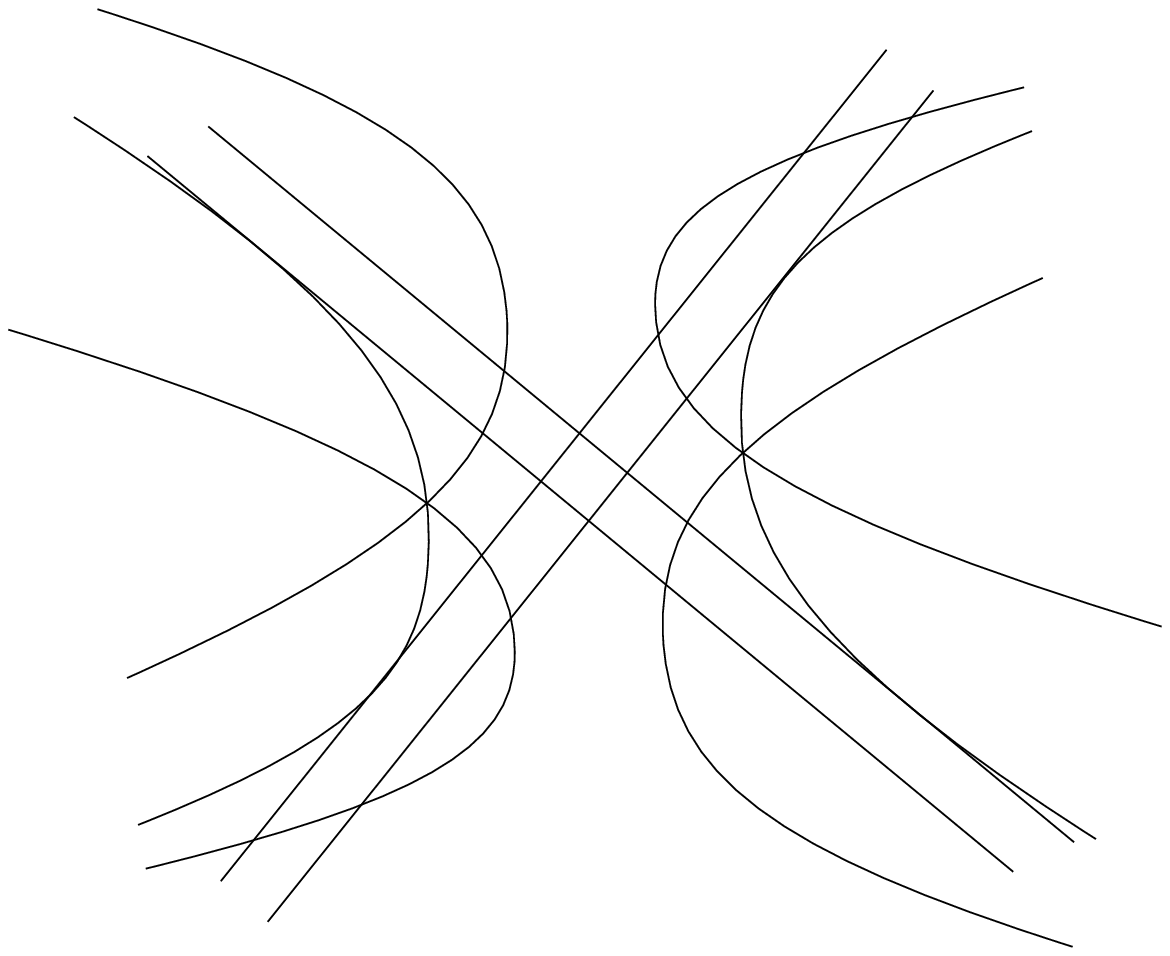}
\caption[]{$S_2$}
\end{center}
\end{figure}


\section{\underline{Computation of the braid monodromy of a curve $S_1$ 
            with}\\
\underline{two conic sections}}\label{S1}
\bthm
Let $S_1$ be a curve over $\R$ in $\C ^2$ of the form:
$\displaystyle S_1 = ({\bigcup _{i=1} ^6} l_i) \cup h_7 \cup ({\bigcup _{i=8} ^9} l_i) \cup h_{10}$,
s.t. $\{ l_i \} _{i=1} ^6,\ \{ l_i \} _{i=8} ^9$ are $8$ lines over $\R$; 
$h_7,h_{10}$ are parabolas over $\R$; $l_4,l_9$ are tangent to the real part 
of $h_{10}$; $l_3,l_8$ are tangent to the real part of $h_7$. 
Let $\{ A_j \} _{j=1} ^{28}$ be singular points of $\pi _1$ as follows:\\
$A_1,A_2 ,A_3,A_5$ are intersection points of the line $l_4$ with the lines
$l_5,l_6,l_8,l_9$ respectively.\\
$A_8,A_9 ,A_{16},A_{23}$ are intersection points of the line $l_3$ with the lines
$l_5,l_6,l_8,l_9$ respectively.\\
$A_{10},A_{15} ,A_{21},A_{26}$ are intersection points of the line $l_2$ 
with the lines $l_5,l_6,l_8,l_9$ respectively.\\
$A_{19},A_{22} ,A_{25},A_{28}$ are intersection points of the line $l_1$ 
with the lines $l_5,l_6,l_8,l_9$ respectively.\\
$A_4,A_6$ are tangent points of the lines $l_9,l_4$ with the parabola $h_{10}$
respectively.\\ 
$A_{12},A_{27}$ are tangent points of the lines $l_3,l_8$ with the parabola $h_7$
respectively.\\
$A_{13},A_{17}$ are intersection points of the lines $l_6,l_5$ with the parabola
$h_7$ respectively.\\
$A_{18},A_{24}$ are intersection points of the line $l_1$ with the parabola
$h_7$.\\
$A_{14},A_{20}$ are intersection points of the line $l_2$ with the parabola
$h_7$.\\
$A_7$ is of the type $a_1$ in the parabola $h_{10}$ and $A_{11}$ is of the type 
$a_2$ in the parabola $h_7$.

Let $N=\{ x(A_j)=x_j \ | \ j=1, \cdots, 28 \}$. Let $E$ be a disk on the x-axis, 
s.t. $N \sbs E - \partial E, \ N \sbs E _{\R}$. Let $M$ be a real point 
on the x-axis, s.t. $x_j \ll M \ \ \forall x_j \in N, \ j=1, \cdots, 28$. 

Then: there exists a g-base $\{ \ell 
(\ga _j) \} _{j=1} ^{28}$ of $\pi _1 (E-N,M)$, 
s.t. each path $\ga _j$ is below the real line and the values 
of $\varphi _M$ w.r.t. this base and $E \times D$ are:

\ethm

\def\uZ{\underline{Z}}
\def\bZ{\bar{Z}}
\def\bZut{\bar{Z}^2}
\def\uZut{\underline{Z}^2}
\def\uZumt{\underline{Z}^{-2}}
\def\bZumt{\bar{Z}^{-2}}

\def\Zoverseven{\mathop{\lower 10pt 
\hbox{$ {\buildrel{\displaystyle\bar{Z}^2_{49}}\over{\hspace{-.2cm}
\scriptstyle{(7)}}} $}}  }

\def\Zoversevenf{\mathop{\lower 10pt 
\hbox{$ {\buildrel{\displaystyle\bar{Z}^2_{4 \; 10}}\over{\hspace{-.5cm}
\scriptstyle{(7)}}} $}}  }

\def\lcZoversevenf{\mathop{\lower 10pt 
\hbox{$ {\buildrel{\displaystyle\bar{z}_{4 \; 10}}\over{\hspace{-.5cm}
\scriptstyle{(7)}}} $}}  }

\def\Zoverst{\mathop{\lower 10pt 
\hbox{$ {\buildrel{\displaystyle\bar{Z}^4_{4 \; 11}}\over{\hspace{-.2cm}
\scriptstyle{(7)(10)}}} $}}  }

\def\Zoverftn{\mathop{\lower 10pt 
\hbox{$ {\buildrel{\displaystyle\bar{Z}^2_{39}}\over{\hspace{-.2cm}
\scriptstyle{(4)}}} $}}  }

\def\lcZoverftn{\mathop{\lower 10pt 
\hbox{$ {\buildrel{\displaystyle\bar{z}_{39}}\over{\hspace{-.2cm}
\scriptstyle{(4)}}} $}}  }

\def\Zovert{\mathop{\lower 10pt 
\hbox{$ {\buildrel{\displaystyle\bar{Z}^2_{2 \; 7}}\over{\hspace{-.4cm}
\scriptstyle{(3)}}} $}}  }

\def\Zovertf{\mathop{\lower 10pt 
\hbox{$ {\buildrel{\displaystyle\bar{Z}^2_{2 \; 9}}\over{\hspace{-.2cm}
\scriptstyle{(3)(4)}}} $}}  }

\def\Zoverftt{\mathop{\lower 10pt 
\hbox{$ {\buildrel{\displaystyle\bar{Z}^2_{3 \; 10}}\over{\hspace{-.5cm}
\scriptstyle{(4)}}} $}}  }

\def\Zovertt{\mathop{\lower 10pt 
\hbox{$ {\buildrel{\displaystyle\bar{Z}^2_{1 \; 7}}\over{\hspace{-.2cm}
\scriptstyle{(2)(3)}}} $}}  }

\def\Zovertfo{\mathop{\lower 10pt 
\hbox{$ {\buildrel{\displaystyle\bar{Z}^2_{1 \; 9}}\over{\hspace{-.2cm}
\scriptstyle{(2)-(4)}}} $}}  }

\def\Zovertft{\mathop{\lower 10pt 
\hbox{$ {\buildrel{\displaystyle\bar{Z}^2_{2 \; 10}}\over{\hspace{-.2cm}
\scriptstyle{(3)(4)}}} $}}  }

\def\lcZovertft{\mathop{\lower 10pt 
\hbox{$ {\buildrel{\displaystyle\bar{z}_{2 \; 10}}\over{\hspace{-.4cm}
\scriptstyle{(3)(4)}}} $}}  }

\def\Zovertfoo{\mathop{\lower 10pt 
\hbox{$ {\buildrel{\displaystyle\bar{Z}^2_{1 \; 10}}\over{\hspace{-.2cm}
\scriptstyle{(2)-(4)}}} $}}  }

\def\Zovere{\mathop{\lower 10pt 
\hbox{$ {\buildrel{\displaystyle\bar{Z}^2_{5 \; 10}}\over{\hspace{-.4cm}
\scriptstyle{(8)}}} $}}  }

\def\lcZovere{\mathop{\lower 10pt 
\hbox{$ {\buildrel{\displaystyle\bar{z}_{5 \; 10}}\over{\hspace{-.4cm}
\scriptstyle{(8)}}} $}}  }

\def\Zoverff{\mathop{\lower 10pt 
\hbox{$ {\buildrel{\displaystyle{Z}^2_{3 \; 6}}\over{\hspace{-.2cm}
\scriptstyle{(5)}}} $}}  }

\def\Zoverfs{\mathop{\lower 10pt 
\hbox{$ {\buildrel{\displaystyle\bar{Z}^2_{2 \; 8}}\over{\hspace{-.2cm}
\scriptstyle{(5)(7)}}} $}}  }

\def\lcZoverfs{\mathop{\lower 10pt 
\hbox{$ {\buildrel{\displaystyle\bar{z}_{2 \; 8}}\over{\hspace{-.2cm}
\scriptstyle{(5)(7)}}} $}}  }

\def\Zoverfive{\mathop{\lower 10pt 
\hbox{$ {\buildrel{\displaystyle\bar{Z}^2_{2 \; 9}}\over{\hspace{-.4cm}
\scriptstyle{(5)}}} $}}  }

\def\Zoverfivetwosix{\mathop{\lower 10pt 
\hbox{$ {\buildrel{\displaystyle\bar{Z}^2_{2 \; 6}}\over{\hspace{-.4cm}
\scriptstyle{(5)}}} $}}  }

\def\zovereight{\mathop{\lower 10pt 
\hbox{$ {\buildrel{\displaystyle\bar{z}_{5 \; 10}}\over{\hspace{-.4cm}
\scriptstyle{(8)}}} $}}  }

\def\Zovereight{\mathop{\lower 10pt 
\hbox{$ {\buildrel{\displaystyle\bar{Z}^2_{5 \; 10}}\over{\hspace{-.4cm}
\scriptstyle{(8)}}} $}}  }

\def\uZ{\underline{Z}}
\def\bZ{\bar{Z}}
\def\bZut{\bar{Z}^2}
\def\uZut{\underline{Z}^2}
\def\uZumt{\underline{Z}^{-2}}
\def\bZumt{\bar{Z}^{-2}}
\noindent
$Z^2_{4\;5} \ , \bZut_{4\;6} \ , \ \Zoverseven \ , \ Z^4_{10 \; 11} \ , \ 
\left(\Zoversevenf \right)^{Z^2_{10 \; 11}} \ , \ \ \Zoverst \ , \ 
\left(Z_{11 \; 11'}\right)^{
\bZut_{4 \; 11} \bZumt_{4 \; 7} Z^2_{10 \; 11}} \ , \ 
\uZut_{3\;5} \ , \ \\
\stackrel{(5)}{Z^2}_{\hspace{-.2cm}3\;6} \ , \ \uZ^2_{2\;5} \ , \ \left(Z_{7 \; 7'} \right)
^{\bZumt_{4\;7} \bZumt_{3\;7}} \ , \ \bZ^4_{3\;7} \ , \uZ^2_{6 \; 7'} \ , \ 
\stackrel{(5)}{\uZ^2}_{ \hspace{-.1cm}2 \; 7'} \ , \ \stackrel{(5)}{\uZ^2}_{ \hspace{-.1cm}2 \; 6} \ , \ 
(\Zoverftn)^{\bZumt_{4\;7}} \ , \ \uZ^2_{5 \; 7'} \ , \uZ^2_{1 \; 7'} \ , \uZ^2_{1 \; 5} \ , \ \\ 
\Zovert \ , \ 
(\Zovertf)^{\bZumt_{4\;7}} \ , \ \stackrel{(5)}{\uZ^2}_{ \hspace{-.1cm}1 \; 6} \ , \left(\Zoverftt\right)^{Z^2_{10 \; 11} \bZumt_{4 \; 7}} \ , \ \Zovertt \ , \ 
\left(\;\;
\Zovertfo \right)^{\bZumt_{4\;7}} \ , \ \left(\;\Zovertft \right)^{\bZumt_{4\;7}Z^2_{10 \; 11}} \ , \\
\\
  \left( \bZ^4_{7 \; 9} \right)^{\bZumt_{4\;7}} \ , \ \left(\;
\Zovertfoo \right)^{\bZumt_{47}  Z^2_{10 \; 11}} .$
\\

\noindent
{\bf Remark:} The paths corresponding to these braids
appear in the proof and in the end of the proof.  The
braids above appear in this way according to section 2.6.

Real($S_1$) is shown in the following figure:
\clearpage
\begin{figure}[h]
\begin{center}
\epsfbox {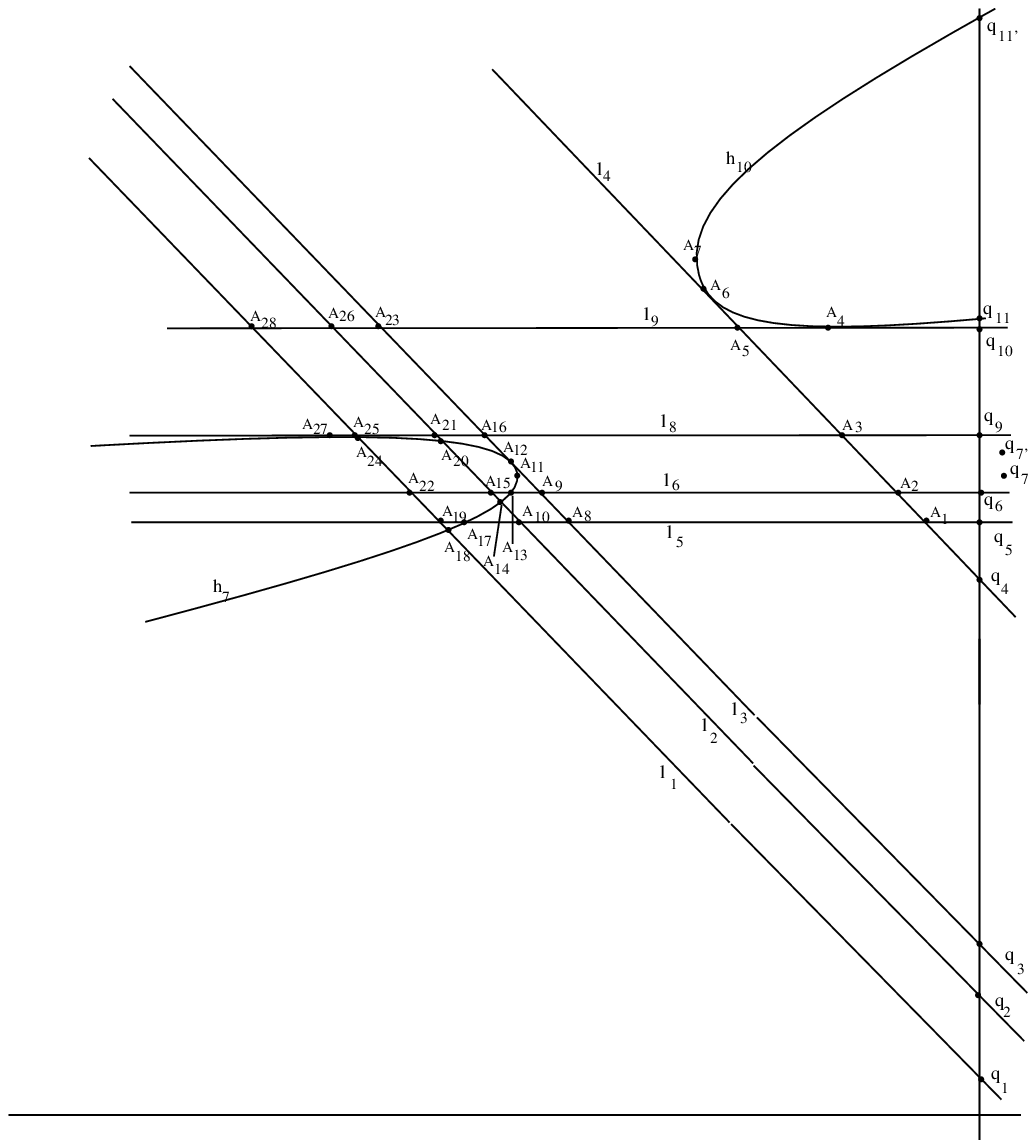}
\caption[]{$S_1$}
\end{center}
\end{figure}


\noindent
{\it Proof:} Let $q_i =l_i \cap K, 1\leq i \leq 6$; $\{ q_7, q_{7'} \} = h_7 \cap K$; 
$q_9 = l_8 \cap K$; $q_{10} = l_9 \cap K$; $\{ q_{11}, q_{11'} \} = h_{10} \cap K$.

Let $K =\{ q_1, \cdots , q_6,q_7, q_{7'}, q_9, q_{10}, q_{11},q_{11'} \}$, s.t. 
$q_1, \cdots , q_6, q_9, q_{10}, q_{11},q_{11'} $ are real and $q_7, q_{7'}$ are complex, 
and:
$q_1< \cdots < q_6< {\rm Re}(q_7)< q_9< q_{10} < q_{11}<q_{11'} 
\ \ \ \ ; \ \ \ \   {\rm Re}(q_7)={\rm Re}(q_{7'}).$
$\be _M$ is defined by: $\be _M (q_i)=i \ \ \ \forall i \leq 6$, 
$\be _M (q_7)=11+i, \ \be _M (q_{7'}) =11-i$ (complex points); $\be _M (q_{11}) =9, \ \be _M (q_{11'}) =10, \ \be _M (q_9) =7, \ \be _M (q_{10}) =8$. Recall that $\pi =\pi _1 | _{S_1} : S_1 \to E$.
Let $\deg \pi = 12, \ \#K_{\R}(x) \geq 8 \ \ \forall x$.

We want to compute $\varphi _M(\ell (\ga _j))$ for $1 \leq j \leq 28$. We choose a
g-base $\{ \ell (\ga _j) \} _{j=1} ^{28}$ of $\pi _1 (E-N , M)$, s.t. each path $\ga _j$
is below the real line. We use the formulas which appear in Theorems 2.25 and 2.28:
${\rm L.V.C.} (\ga _j)= (\xi _{x_j '}) \Psi _{\ga _j '}\ \ , \qquad 
\varphi _M(\ell (\ga _j))= \De < {\rm L.V.C.}(\ga _j)> ^{\ep _{x_j}}.$

Using these formulas, we need $\lambda_{x_j}, \epsilon_{x_j}, \delta_{x_j}$ for
$1 \leq j \leq 28$. Having the following table, we can apply the formulas above.
\begin{center}
\begin{tabular}{|c|c|c|c|}
\hline
j & $\la _{x_j}$ & $\ep _{x_j}$ & $\de _{x_j}$ \\
\hline
1 & $<4,5>$ & 2 & $\De <4,5>$ \\
2 & $<5,6>$ & 2 & $\De <5,6>$ \\
3 & $<6,7>$ & 2 & $\De <6,7>$ \\
4 & $<8,9>$ & 4 & $\De^2 <8,9>$ \\
5 & $<7,8>$ & 2 & $\De <7,8>$ \\
6 & $<8,9>$ & 4 & $\De^2 <8,9>$ \\
7 & $<9,10>$ & 1 & $\De^{\frac{1}{2}'}_{I_4 I_2}  <9>$ \\
8 & $<3,4>$ & 2 & $\De <3,4>$ \\
9 & $<4,5>$ & 2 & $\De <4,5>$ \\
10 & $<2,3>$ & 2 & $\De <2,3>$ \\
11 & $P_5$ & 1  &  $\De _{I_2 I_4} ^{1 \over 2} <5>$ \\
12 & $<6,7>$ & 4 & $\De^2 <6,7>$ \\
13 & $<4,5>$ & 2 & $\De <4,5>$ \\
14 & $<3,4>$ & 2 & $\De <3,4>$ \\
15 & $<4,5>$ & 2 & $\De <4,5>$ \\
16 & $<7,8>$ & 2 & $\De <7,8>$ \\
17 & $<2,3>$ & 2 & $\De <2,3>$ \\
18 & $<1,2>$ & 2 & $\De <1,2>$ \\
19 & $<2,3>$ & 2 & $\De <2,3>$ \\
20 & $<5,6>$ & 2 & $\De <5,6>$ \\
21 & $<6,7>$ & 2 & $\De <6,7>$ \\
22 & $<3,4>$ & 2 & $\De <3,4>$ \\
23 & $<8,9>$ & 2 & $\De <8,9> $\\
24 & $<4,5>$ & 2 & $\De <4,5>$ \\
25 & $<5,6>$ & 2 & $\De <5,6>$ \\
26 & $<7,8>$ & 2 & $\De <7,8> $\\
27 & $<4,5>$ & 4 & $\De ^2<4,5>$ \\
28 & $<6,7>$ & 2 & $\De <6,7>$ \\
\hline
\end{tabular}
\end{center}

We do not present $28$ computations, but just the $4$ ones for
$j = 5,11,16,26$. For all the others, we present just the final results of the 
computations. For all these missing computations, see [Am].

\medskip
\noindent
$(\xi_{x'_5}) \Psi_{\ga'_5}   = 
<7,8>\De^2<8,9>\De<6,7>\De<5,6>\De<4,5> \beta^{-1}_M = (\lcZoversevenf)^{Z^2_{10 \; 11}}$

\medskip
\noindent
$<7,8>$ $\hspace{1cm}\vcenter{\hbox{\epsfbox{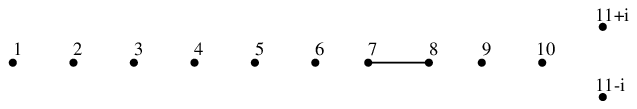}}}$

\medskip
\noindent
$\De^2<8,9>$ $\hspace{.4cm}\vcenter{\hbox{\epsfbox{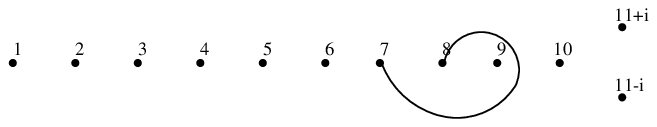}}}$

\medskip
\noindent
$\De<6,7>$ $\hspace{.5cm}\vcenter{\hbox{\epsfbox{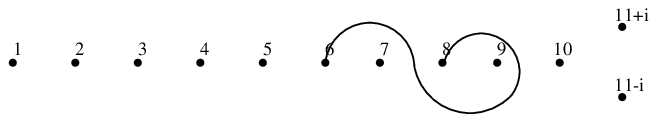}}}$

\medskip
\noindent
$\De<5,6>$ $\hspace{.5cm}\vcenter{\hbox{\epsfbox{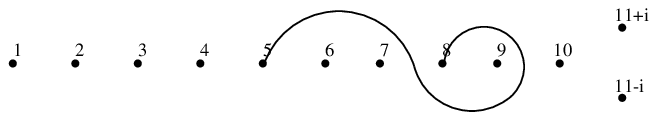}}}$

\medskip\noindent
$\De<4,5>$ $\hspace{.5cm}\vcenter{\hbox{\epsfbox{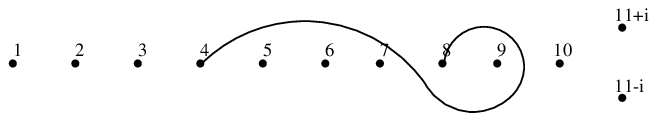}}}$

\medskip
\noindent
$\beta^{-1}_M$ $\hspace{1.7cm}\vcenter{\hbox{\epsfbox{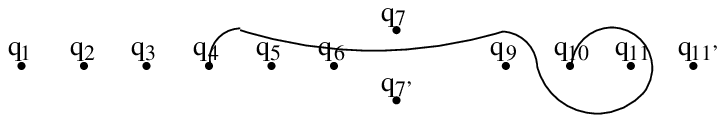}}}$

\medskip
\noindent
$\vp_M(\ell(\ga_5)) = (\Zoversevenf)^{Z^2_{10 \; 11}}$\\

\medskip
\noindent
$(\xi_{x'_{11}}) \Psi_{\ga'_{11}}   = P_5 \De
<2,3>\De<4,5>\De<3,4>\De^{\frac{1}{2}'}_{I_4I_2} <9> \De^2<8,9> 
\De<7,8>\De^2<8,9>\De<6,7>\De<5,6>\De<4,5> \beta^{-1}_M = \left(z_{7 \; 7'} \right)
^{\bZumt_{4\;7} \bZumt_{3\;7}}$

\medskip\noindent
$P_5$ $\hspace{1.8cm}\vcenter{\hbox{\epsfbox{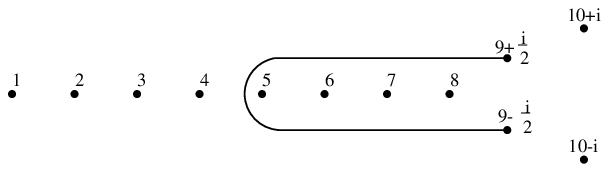}}}$

\medskip\noindent
$\De<2,3>$ doesn't change it

\medskip\noindent
$\De<4,5>$ $\hspace{.5cm}\vcenter{\hbox{\epsfbox{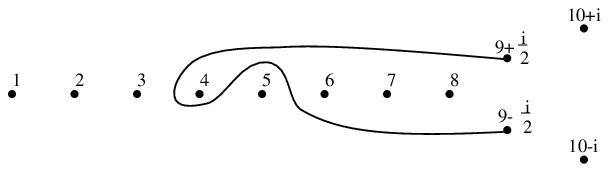}}}$

\medskip\noindent
$\De<3,4>$ $\hspace{.5cm}\vcenter{\hbox{\epsfbox{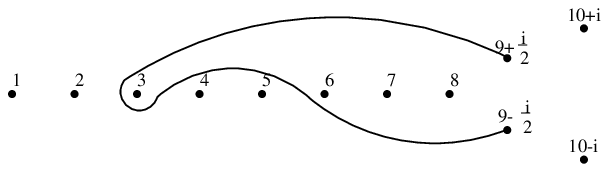}}}$

\medskip\noindent
$\De ^{\frac{1}{2}'}_{I_4I_2} <9>$
$\hspace{.5cm}\vcenter{\hbox{\epsfbox{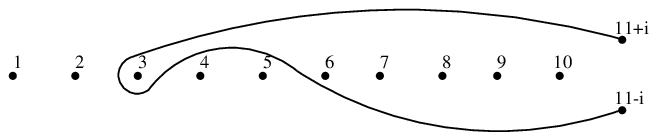}}}$

\medskip\noindent
$ \De^2<8,9> 
\De<7,8>\De^2<8,9>\De<6,7>$ don't change it

\medskip
\noindent
$\De<5,6>$ $\hspace{.5cm}\vcenter{\hbox{\epsfbox{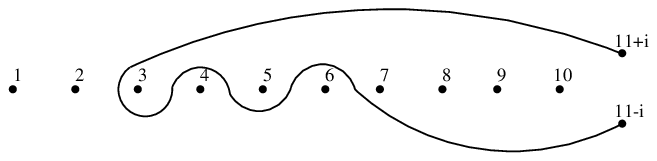}}}$

\medskip\noindent
$\De<4,5>$ $\hspace{.5cm}\vcenter{\hbox{\epsfbox{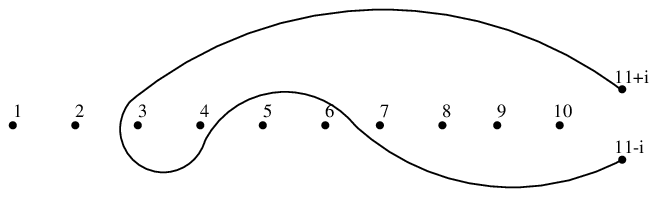}}}$

\medskip\noindent
$\beta^{-1}_M$ $\hspace{2cm}\vcenter{\hbox{\epsfbox{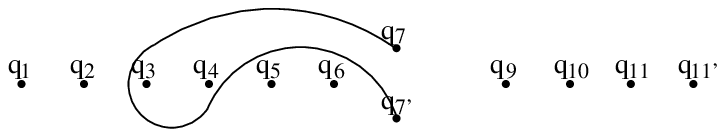}}}$

\medskip\medskip
\noindent
$\vp_M(\ell(\ga_{11}))=
\left(Z_{7 \; 7'} \right)
^{\bZumt_{4\;7} \bZumt_{3\;7}}$\\

\medskip
\noindent
$(\xi_{x'_{16}}) \Psi_{\ga'_{16}}   = <7,8>\De<4,5>\De<3,4>
\De<4,5>\De^2<6,7> \De^{\frac{1}{2}}_{I_2I_4} <5> \De<2,3>\\
\De<4,5>\De<3,4> \De^{\frac{1}{2}'}_{I_4I_2} <9>  \De^2<8,9> 
\De<7,8>\De^2<8,9>\De<6,7>\De<5,6>\\
\De<4,5> \beta^{-1}_M = \left
(\lcZoverftn\right)^{\hspace{-.1cm}\bZumt_{4\;7}}$

\medskip
\noindent
$<7,8>$ $\hspace{1cm}\vcenter{\hbox{\epsfbox{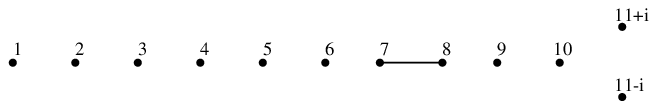}}}$

\medskip\noindent
$\De<4,5>\De<3,4>
\De<4,5>$ don't change it\\

\medskip
\noindent
$\De^2<6,7>$ $\hspace{1cm}\vcenter{\hbox{\epsfbox{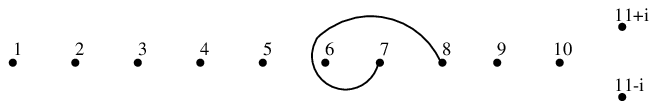}}}$

\medskip\noindent
$\De^{\frac{1}{2}}_{I_2I_4} <5>$ 
$\hspace{1cm}\vcenter{\hbox{\epsfbox{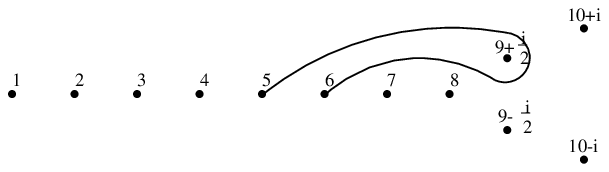}}}$

\medskip
\noindent
$\De<2,3> $ doesn't change it\\

\medskip
\noindent
$\De<4,5>$ $\hspace{1cm}\vcenter{\hbox{\epsfbox{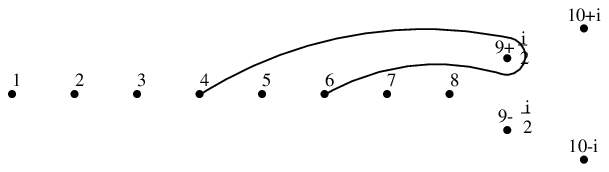}}}$

\medskip\noindent
$\De<3,4>$ $\hspace{1cm}\vcenter{\hbox{\epsfbox{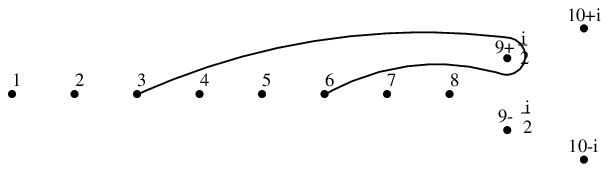}}}$

\noindent
$ \De^{\frac{1}{2}'}_{I_4I_2} <9>$
$\hspace{1cm}\vcenter{\hbox{\epsfbox{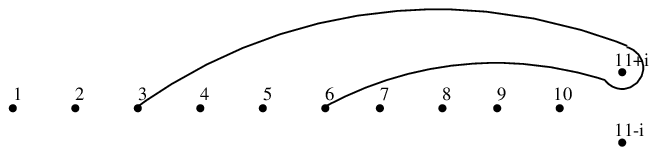}}}$

\noindent
$ \De^2<8,9> 
\De<7,8>\De^2<8,9>$  don't change it\\

\medskip
\noindent
$\De<6,7>$ $\hspace{1cm}\vcenter{\hbox{\epsfbox{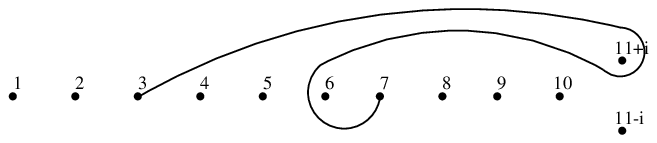}}}$

\noindent
$\De<5,6> $ $\hspace{1cm}\vcenter{\hbox{\epsfbox{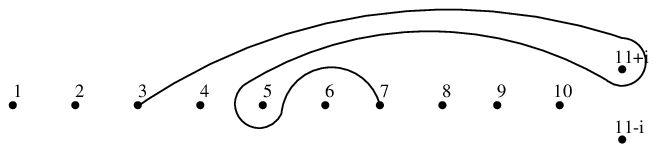}}}$

\noindent
$\De<4,5>$ $\hspace{1cm}\vcenter{\hbox{\epsfbox{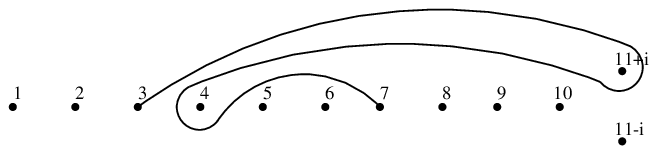}}}$

\noindent
$\beta^{-1}_M$ $\hspace{1.7cm}\vcenter{\hbox{\epsfbox{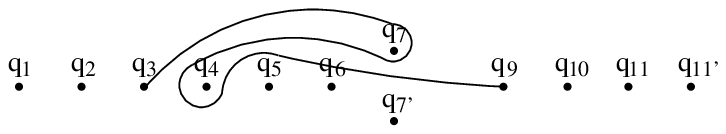}}}$

\medskip
\noindent
$\vp_M(\ell(\ga_{16}))= (\Zoverftn)^{\bZumt_{4\;7}}$\\

\medskip
\noindent
$(\xi_{x'_{26}}) \Psi_{\ga'_{26}}   = <7,8>\De<5,6>\De<4,5>\De
<8,9>\De<3,4>\De<6,7>\De<5,6>\\
\De<2,3>\De<1,2>\De<2,3>\De
<7,8>\De<4,5>\De<3,4>
\De<4,5>\De^2<6,7> \De^{\frac{1}{2}}_{I_2I_4} <5> \De<2,3>
\De<4,5>\De<3,4> \De^{\frac{1}{2}'}_{I_4I_2} <9>  \De^2<8,9> 
\De<7,8>\De^2<8,9>\De<6,7>\De<5,6>\De<4,5> \beta^{-1}_M = \left
(\;\;\lcZovertft\right)^{\bZumt_{4\;7}Z^2_{10 \; 11}}$\\

\medskip
\noindent
$<7,8>$ $\hspace{.5cm}\vcenter{\hbox{\epsfbox{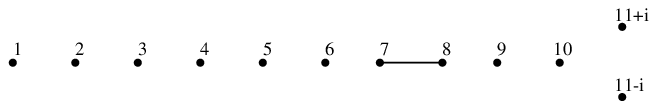}}}$

\medskip
\noindent
$\De<5,6>\De<4,5>$ don't change it\\

\medskip
\noindent
$<8,9>$ $\hspace{.5cm}\vcenter{\hbox{\epsfbox{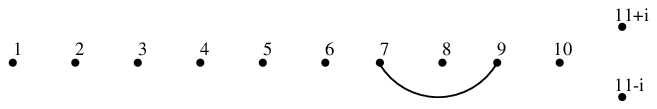}}}$

\medskip
\noindent
$\De<3,4> $ doesn't change it\\

\noindent
$\De<6,7>$ $\hspace{.5cm}\vcenter{\hbox{\epsfbox{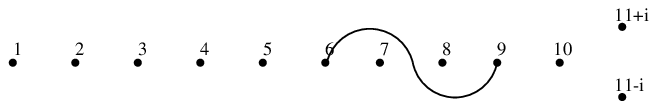}}}$

\medskip
\noindent
$ \De<5,6>$ $\hspace{.5cm}\vcenter{\hbox{\epsfbox{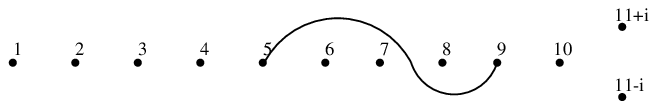}}}$

\medskip
\noindent
$\De<2,3>\De<1,2>\De<2,3> $ don't change it\\

\medskip
\noindent
$\De<7,8>$ $\hspace{.5cm}\vcenter{\hbox{\epsfbox{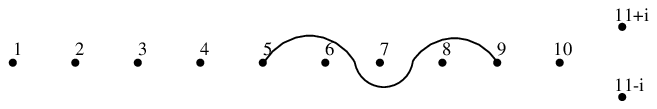}}}$

\medskip
\noindent
$\De<4,5>$ $\hspace{.5cm}\vcenter{\hbox{\epsfbox{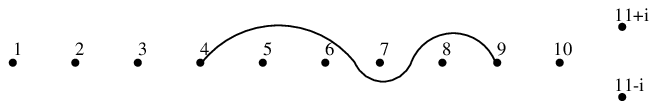}}}$

\medskip
\noindent
$\De<3,4>$ $\hspace{.5cm}\vcenter{\hbox{\epsfbox{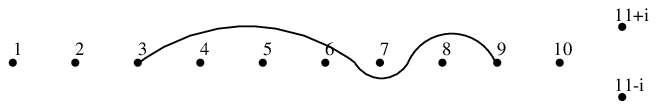}}}$

\medskip
\noindent
$\De<4,5>$ doesn't change it\\

\medskip
\noindent
$\De^2<6,7>$ $\hspace{.5cm}\vcenter{\hbox{\epsfbox{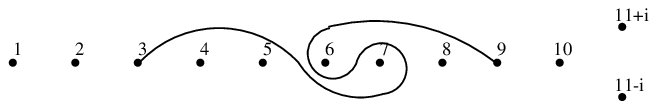}}}$

\medskip
\noindent
$ \De^{\frac{1}{2}}_{I_2I_4} <5>$
$\hspace{.5cm}\vcenter{\hbox{\epsfbox{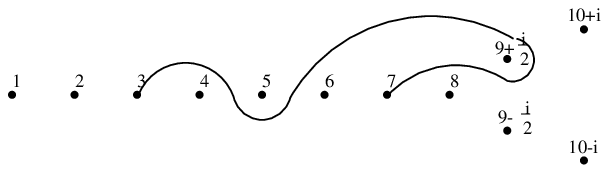}}}$

\noindent
$\De<2,3>$ $\hspace{.5cm}\vcenter{\hbox{\epsfbox{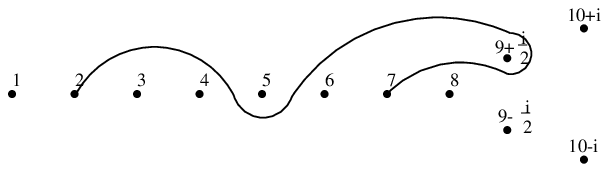}}}$

\noindent
$\De<4,5>$ $\hspace{.5cm}\vcenter{\hbox{\epsfbox{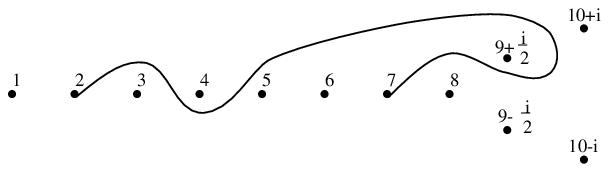}}}$

\noindent
$\De<3,4>$ $\hspace{.5cm}\vcenter{\hbox{\epsfbox{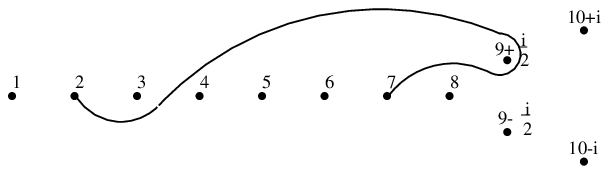}}}$

\noindent
$\De^{\frac{1}{2}'}_{I_4I_2} <9>  $
$\hspace{.5cm}\vcenter{\hbox{\epsfbox{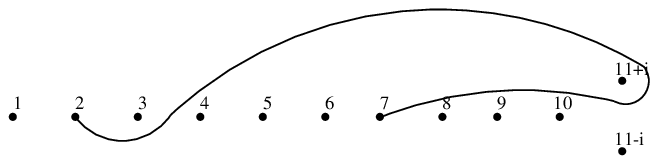}}}$

\noindent
$\De^2<8,9> $ doesn't change it\\

\noindent
$\De<7,8>$ $\hspace{.5cm}\vcenter{\hbox{\epsfbox{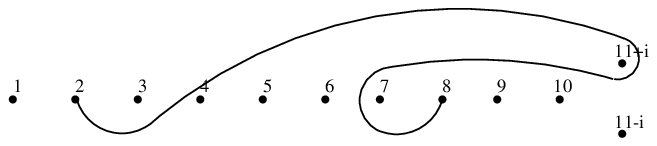}}}$

\noindent
$\De^2<8,9>$ $\hspace{.5cm}\vcenter{\hbox{\epsfbox{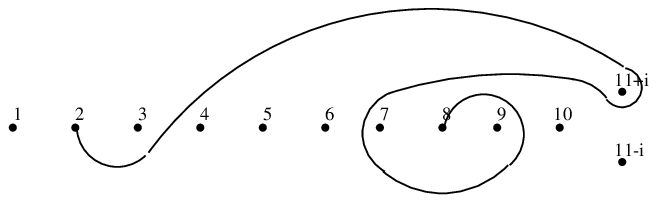}}}$

\noindent
$\De<6,7>$ $\hspace{.5cm}\vcenter{\hbox{\epsfbox{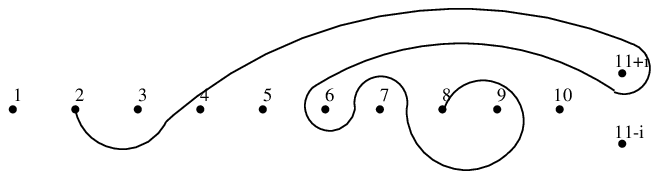}}}$

\noindent
$\De<5,6>$ $\hspace{.5cm}\vcenter{\hbox{\epsfbox{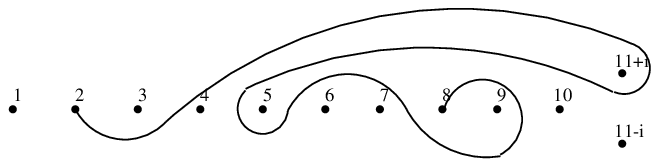}}}$

\noindent
$\De<4,5>$ $\hspace{.5cm}\vcenter{\hbox{\epsfbox{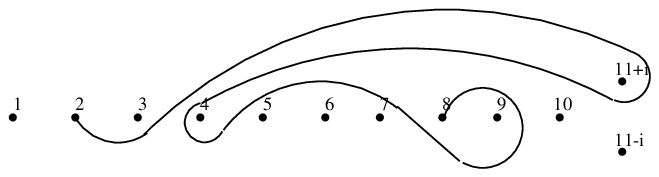}}}$

\noindent
$\beta^{-1}_M$ $\hspace{1.7cm}\vcenter{\hbox{\epsfbox{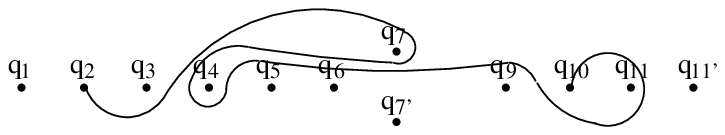}}}$

\noindent
$\vp_M(\ell(\ga_{26}))= \left
(\Zovertft\right)^{\bZumt_{4\;7}Z^2_{10 \; 11}}$\\

Here are the final results of the other computations: \\

\medskip

\noindent
$\vcenter{\hbox{\epsfbox{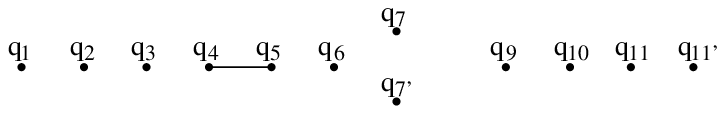}}}$\hspace{3cm}
$\vp_M(\ell(\ga_{1})) = Z^2_{4\;5}$

\medskip
\noindent
$\vcenter{\hbox{\epsfbox{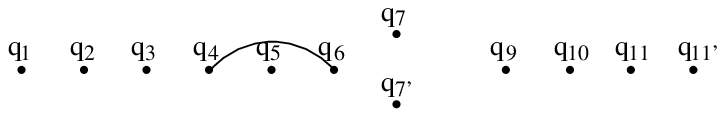}}}$\hspace{3cm}
$\vp_M(\ell(\ga_{2})) = \bZut_{4\;6}$

\medskip

\noindent
$\vcenter{\hbox{\epsfbox{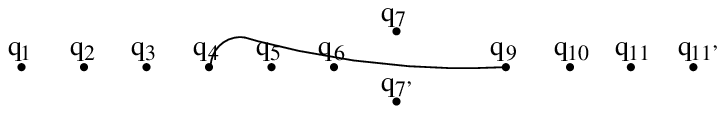}}}$\hspace{3cm}
$\vp_M(\ell(\ga_{3})) = \Zoverseven$

\medskip

\noindent
$\vcenter{\hbox{\epsfbox{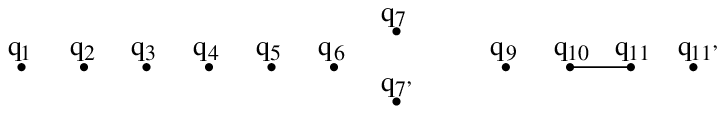}}}$\hspace{3cm}
 $\vp_M(\ell(\ga_{4})) = Z^4_{10 \; 11} $

\medskip

\noindent
$\vcenter{\hbox{\epsfbox{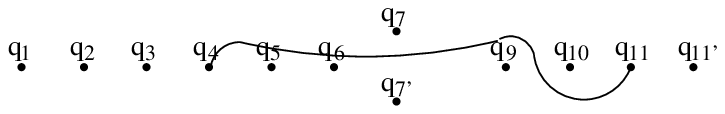}}}$\hspace{3cm}
$\vp_M(\ell(\ga_{6})) = \Zoverst  $

\medskip

\noindent
$\vcenter{\hbox{\epsfbox{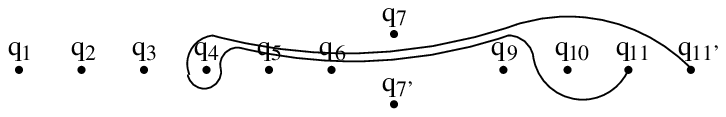}}}$\hspace{3cm}
$\vp_M(\ell(\ga_{7})) = \left(Z_{11 \; 11'}\right)^{
\bZut_{4 \; 11} \bZumt_{4 \; 7} Z^2_{10 \; 11}}$

\medskip

\noindent
$\vcenter{\hbox{\epsfbox{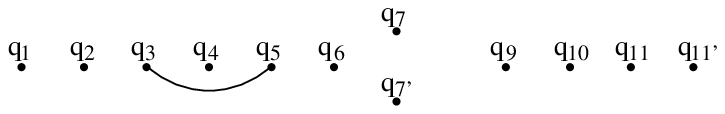}}}$\hspace{3cm} 
$ \vp_M(\ell(\ga_{8})) = \uZ^2_{3 \; 5}$

\medskip

\noindent
$\vcenter{\hbox{\epsfbox{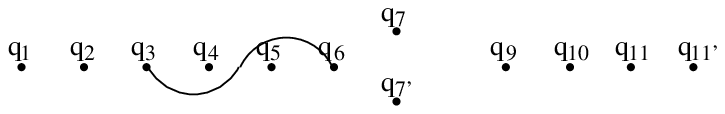}}}$\hspace{3cm} 
$ \vp_M(\ell(\ga_{9})) = \stackrel{(5)}{Z^2}_{\hspace{-.2cm}3\;6} $

\medskip

\noindent
$\vcenter{\hbox{\epsfbox{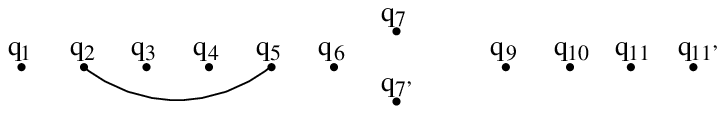}}}$\hspace{3cm} 
$ \vp_M(\ell(\ga_{10})) = \uZ^2_{2\;5} $

\medskip
\noindent
$\vcenter{\hbox{\epsfbox{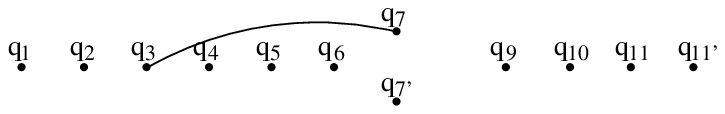}}}$\hspace{3cm} 
$ \vp_M(\ell(\ga_{12})) = \bZ^4_{3\;7} $

\medskip

\noindent
$\vcenter{\hbox{\epsfbox{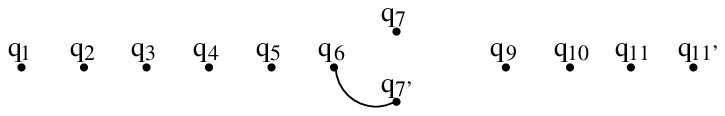}}}$\hspace{3cm} 
$\vp_M(\ell(\ga_{13}))= \uZ^2_{6 \; 7'}$

\medskip

\noindent
$\vcenter{\hbox{\epsfbox{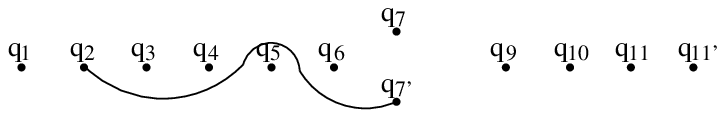}}}$\hspace{3cm} 
$ \vp_M(\ell(\ga_{14})) = \stackrel{(5)}{\uZ^2}_{ \hspace{-.1cm}2 \; 7'}$

\medskip

\noindent
$\vcenter{\hbox{\epsfbox{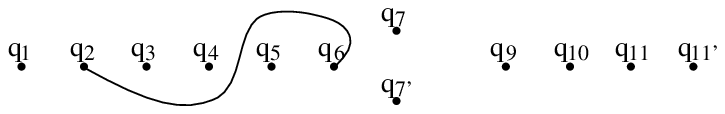}}}$\hspace{3cm} 
$ \vp_M(\ell(\ga_{15})) = \stackrel{(5)}{\uZ^2}_{ \hspace{-.1cm}2 \; 6}$

\medskip

\noindent
$\vcenter{\hbox{\epsfbox{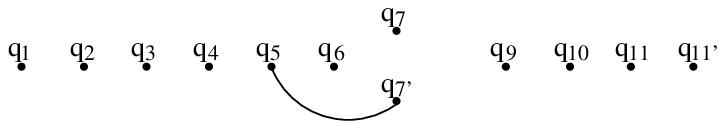}}}$\hspace{3cm} 
$ \vp_M(\ell(\ga_{17}) ) = \uZ^2_{5 \; 7'}$

\medskip

\noindent
$\vcenter{\hbox{\epsfbox{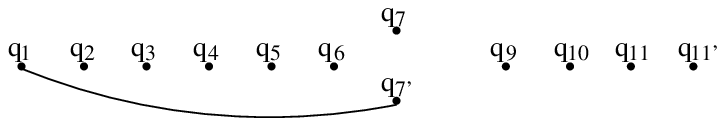}}}$\hspace{3cm} 
$ \vp_M(\ell(\ga_{18}) ) = \uZ^2_{1 \; 7'} $

\medskip

\noindent
$\vcenter{\hbox{\epsfbox{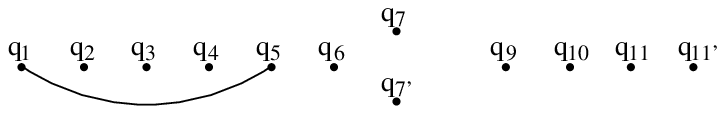}}}$\hspace{3cm} 
$ \vp_M(\ell(\ga_{19}) ) =  \uZ^2_{1 \; 5}$

\medskip

\noindent
$\vcenter{\hbox{\epsfbox{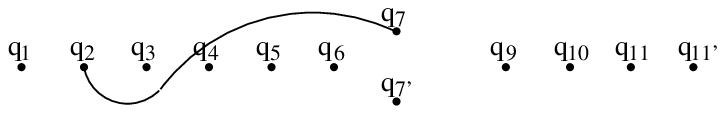}}}$\hspace{3cm} 
$\vp_M(\ell(\ga_{20})) =  \Zovert$

\medskip

\noindent
$\vcenter{\hbox{\epsfbox{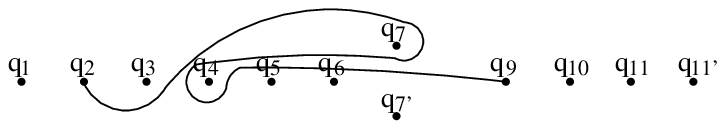}}}$\hspace{3cm} 
$ \vp_M(\ell(\ga_{21})) =  \left(\;\;\Zovertf
\right)^{\bZumt_{4\;7}}$

\medskip

\noindent
$\vcenter{\hbox{\epsfbox{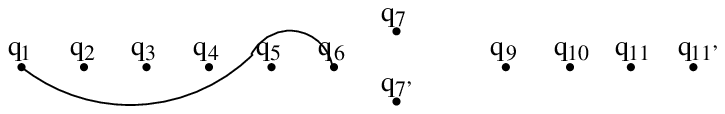}}}$\hspace{3cm} 
$\vp_M(\ell(\ga_{22}) ) = \stackrel{(5)}{\uZ^2}_{ \hspace{-.1cm}1 \; 6}$

\medskip

\noindent
$\vcenter{\hbox{\epsfbox{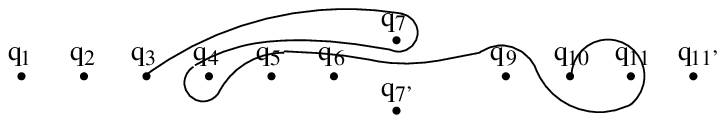}}}$\hspace{3cm} 
$ \vp_M(\ell(\ga_{23}) ) = 
\left(\Zoverftt \right)^{Z^2_{10 \; 11} \bZumt_{4 \; 7}}$

\medskip

\noindent
$\vcenter{\hbox{\epsfbox{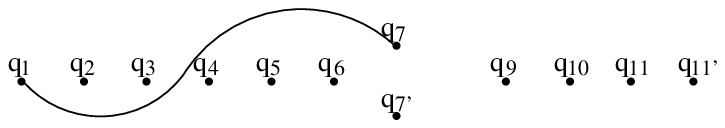}}}$\hspace{3cm} 
$ \vp_M(\ell(\ga_{24})) = \Zovertt$

\medskip

\noindent
$\vcenter{\hbox{\epsfbox{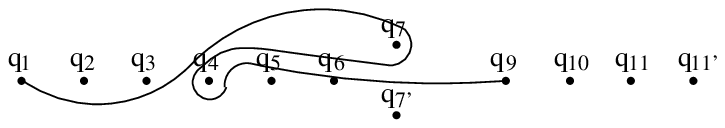}}}$\hspace{3cm} 
$ \vp_M(\ell(\ga_{25})) = \left(\;\;\Zovertfo\right)^{\bZumt_{4\;7}}$

\medskip

\noindent
$\vcenter{\hbox{\epsfbox{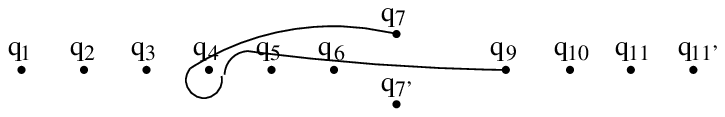}}}$\hspace{3cm} 
$ \vp_M(\ell(\ga_{27}) ) = \left( \bZ^4_{7 \; 9} \right)^{\bZumt_{4\;7}}$

\medskip

\noindent
$\vcenter{\hbox{\epsfbox{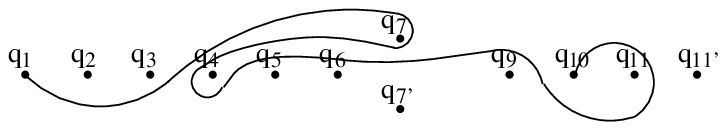}}}$\hspace{3cm} 
$ \vp_M(\ell(\ga_{28}) ) = \left(\;\;\Zovertfoo \right)^{\bZumt_{47} Z^2_{10 \; 11}}$\\

\vspace{-1cm}

\begin{flushright}
$\Box$
\end{flushright}

\section{\underline{Computation of the braid monodromy of a curve $S_2$ with}\\
\underline{three conic sections}}\label{S2}

\bthm
Let $S_2$ be a curve over $\R$ in $\C ^2$ of the form:
$S_2 = h_1 \cup l_2 \cup  l_3 \cup  h_4 \cup h_5 \cup  l_6
\cup  l_7$,
s.t. $ l_2,l_3,l_6,l_7 $ are lines over $\R$; $ h_1,h_4,h_5 $ 
are hyperbolas over $\R$; the lines $l_2,l_3,l_6,l_7$ are tangent to 
the real part of $h_4$;
the lines $l_2,l_3$ intersect in a point $1 \ll x(l_2 \cap l_3);$ 
the lines  $l_6,l_7$ intersect in a point $x(l_6 \cap l_7)$, s.t. 
$x (l_2 \cap l_3) <x(l_6 \cap l_7)$; $x (l_2 \cap l_3), x(l_6 \cap l_7) \not \in E$ for $E$ 
a disk on the x-axis.

Let $\{ A_j \} _{j=1} ^{28}$ be singular points of $\pi _1$ as follows: \\
$A_1,A_4$ are tangent points of the lines $l_3,l_6$ with the hyperbola $h_4$ 
respectively.\\ 
$A_2,A_7$ are intersection points of the hyperbola $h_5$ with the line $l_6$.\\
$A_3,A_{11}$ are intersection points of the hyperbola $h_5$ with the line $l_7$.\\
$A_8,A_9$ are intersection points of the lines $l_3,l_2$ with the hyperbola 
$h_1$ respectively.\\
$A_{13},A_{14}$ are intersection points of the lines $l_3,l_2$ with the line $l_6$ 
 respectively.\\
$A_{15},A_{16}$ are intersection points of the lines $l_3,l_2$ with the line $l_7$ 
 respectively.\\
$A_{18},A_{26}$ are intersection points of the hyperbola $h_5$ with the line $l_6$. \\
$A_{20},A_{22}$ are intersection points of the lines $l_3,l_2$ with the 
hyperbola $h_1$ respectively. \\ 
$A_{25},A_{28}$ are tangent points of the lines $l_7,l_2$ with the hyperbola $h_4$ 
 respectively.\\
$A_{21},A_{27}$ are intersection points of the hyperbola $h_5$ with 
the line $l_7$.\\ 
$A_5,A_{24}$ are intersection points of the hyperbolas $h_1,h_4,h_5$.\\
$A_6,A_{10},A_{12}$ are points of the type $a_1$ of the hyperbolas $h_4,h_1,h_5$ 
respectively.\\
$A_{17},A_{19},A_{23}$ are points of the type $a_2$ of the hyperbolas $h_5,h_1,h_4$ 
respectively.

Let $N=\{ x(A_j)=x_j \ | \ 1 \leq j \leq 28 \}$. Let $E$ be a disk on the x-axis, 
s.t. $N \sbs E- \partial E, N \sbs E _{\R}$.

Let $M$ be a real point on the x-axis, s.t. 
$x_j \ll M  \ \ \forall x_j \in N, 1 \leq j \leq 28$.

Then: there exists a g-base $\{ \ell (\ga _j) \} _{j=1} ^{28}$ of $\pi _1 (E-N,M)$, 
s.t. the first 14 paths $\ga _j$ are below the real line and the last 14 
paths $\ga _j, 15 \leq j \leq 28$, have a part below the real line and a part above 
the real line (see Figure \ref{gbase}).

Moreover, $\varphi _M$, the braid monodromy w.r.t. $E \times D$, 
is $\varphi _M = F_1 \cdot F_1 ^{\rho^{-1}}$, where:

\medskip
\noindent
$F_1=   Z^4_{3 \; 4} \cdot Z^2_{8 \; 9} \cdot \bZ^2_{8 \; 10}    \cdot 
\uZ^4_{7 \; 9} \cdot \left((Z^2_{45})^{Z^2_{34}} \ , Z^2_{5 \; 6}\right) 
\cdot \left(\bZ_{4 \; 7} \right)^{Z^2_{34}\uZ^2_{79}} \cdot 
\stackrel{(6)}{\uZ^2}_{\hspace{-.2cm}5 \; 9} \cdot 
\stackrel{(4)}{Z^2}_{\hspace{-.2cm}3 \; 6}
 \cdot \uZ^2_{2 \; 6} \cdot \uZ_{1 \; 6} \cdot 
\Zovere
\cdot \bZ_{5 \; 8} \cdot \stackrel{(4)(6)}{\uZ^2}_{\hspace{-.3cm}3 \; 9} \cdot
\stackrel{(6)}{\uZ^2}_{\hspace{-.1cm}2 \; 9}$,
s.t. $\rho = (\De <2,3> \De <9,10>)$.
\ethm
Real$(S_2)$ is shown in the following figure:

\vspace{-1cm}
\begin{figure}[htp]
\begin{center}
\setlength{\epsfysize}{8.7cm}
\epsfbox {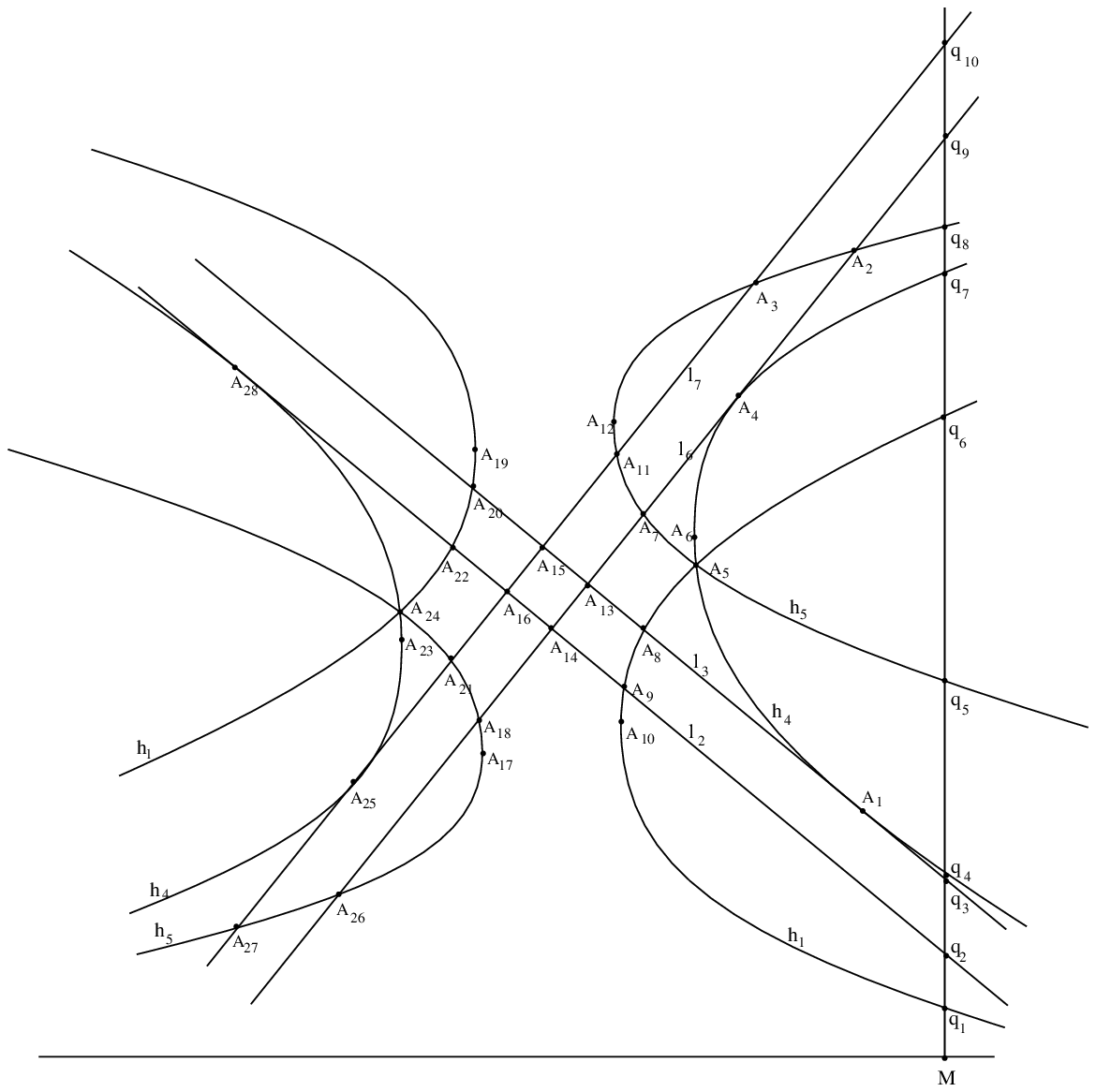}
\caption[]{$S_2$} \label{s2fig}
\end{center}
\end{figure}


\noindent
{\it Proof:} Let us choose a g-base of $\pi _1 (E-N,M)$ as follows: let us take 
the first 14 elements of the g-base the loops constructed from the standard
bush from $M$ to $\{ x_j \} _{j=1} ^{14}$ below the real line, namely:
$\ell(\ga _1), \cdots , \ell (\ga _{14})$.

\vspace{-0.1cm}
For the last $14$ elements of the g-base we take 
$\ell (\ga _{15}), \cdots , \ell (\ga _{28})$ where $\ga _j$ are constructed as follows: let $-P,P_1, P$ be points on the real line which satisfy
$M \ll x(l_2 \cap l_3) < P_1 < x(l_6 \cap l_7) < P \ \ \ , \qquad -P \ll x_{28}$. 
Let $T$ be a big semicircle below the real line from $-P$ to $P$; $T_1$ 
be a semicircle below the real line from $P$ to $P_1$; $T_2$ be 
a semicircle below the real line from $P_1$ to $M$.

Let $\tilde \ga _j$ be paths above the real line from $x_j$ to $-P$, 
$15 \leq j \leq 28$. Let $\ga _j$ be the paths of the form: 
$\ga _j = \tilde \ga _j T T_1 T_2, \ 15 \leq j \leq 28$.

\begin{figure}[htp]
\begin{center}
\epsfbox {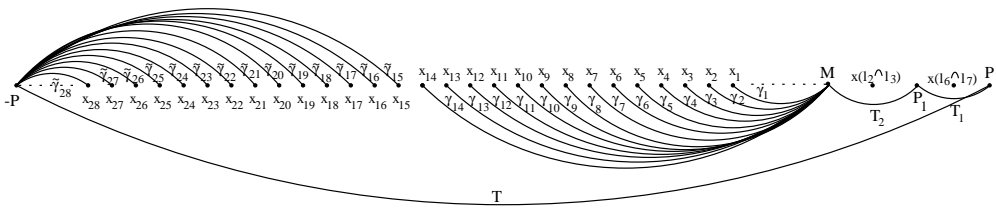}
\caption[]{The g-base} \label{gbase}
\end{center}

\vspace{-1cm}

\end{figure}


We use the formulas:
${\rm L.V.C.} (\ga _j)= (\xi _{x_j '}) \Psi _{\ga _j '}\ , \
\varphi _M(\ell (\ga _j))= \De < {\rm L.V.C.}(\ga _j)> 
^{\ep _{x_j}}$.  In order to use them we need to compute 
$\la _{x_j},\ep _{x_j} ,\de _{x_j}$ for $1 \leq j \leq 14$.
\vspace{-0.2cm}
We present all of them in the following table:

\begin{center}
\begin{tabular}{|c|c|c|c|}
\hline
j & $\la _{x_j}$ & $\ep _{x_j}$ & $\de _{x_j}$ \\
\hline
1 & $<3,4>$ & 4 & $\De ^2<3,4>$ \\
2 & $<8,9>$ & 2 & $\De <8,9>$ \\
3 & $<9,10>$ & 2 & $\De <9,10>$ \\
4 & $<7,8>$ & 4 & $\De ^2 <7,8>$ \\
5 & $<4,6>$ & 2 & $\De <4,6>$ \\
6 & $<6,7>$ & 1 & $\De ^{1 \over 2} _{I_2 \R} <6>$ \\
7 & $<5,6>$ & 2 & $\De <5,6>$ \\
8 & $<3,4>$ & 2 & $\De <3,4>$ \\
9 & $<2,3>$ & 2 & $\De <2,3>$\\
10 & $<1,2>$ & 1 & $\De ^{1\over 2} _{I_4 I_2} <1>$\\
11 & $<4,5>$ & 2 & $\De <4,5>$\\
12 & $<5,6>$ & 1 & $\De ^{1 \over 2} _{I_6 I_4} <5>$\\
13 & $<2,3>$ & 2 & $\De <2,3>$\\
14 & $<1,2>$ & 2 & $\De <1,2>$\\
\hline
\end{tabular}
\end{center}
\clearpage

By Figure \ref{s2fig}, 
$h_1 \cap K = \{ q_1, q_6 \}, \ h_4 \cap K = \{ q_4, q_7 \}, \ h_5 \cap K = \{ q_5, q_8 \}, \ l_7 \cap K = q_{10} , \ l_6 \cap K = q_9, \ l_3 \cap K = q_3, \ l_2 \cap K = q_2$. 

So, $K=\{ q_1, \cdots , q_{10} \}, \ q_1 < \cdots < q_{10}$. We choose a diffeomorphism 
$\be _M$ which satisfies: $\be _M (q_i) = i$.

We do not present here all the $14$ computations but a few ones,
 for 
$j=5,11,13$. For the others, we present just the final results. For all 
the details, see [Am].
\\
\\
\noindent
$(\xi_{x'_5}) \Psi_{\ga'_5}   = 
<4,6>\De^2<7,8>\De<9,10>\De<8,9>\De^2<3,4> \beta^{-1}_M = 
((z_{45})^{Z^2_{34}} , z_{56})$

\medskip
\noindent
$<4,6>$ $\hspace{1cm}\vcenter{\hbox{\epsfbox{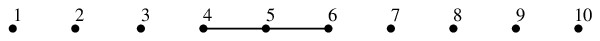}}}$

\medskip
\noindent
$\De^2<7,8>\De<9,10>\De<8,9>$ don't change it

\medskip
\noindent
$\De^2<3,4>$ $\hspace{1cm}\vcenter{\hbox{\epsfbox{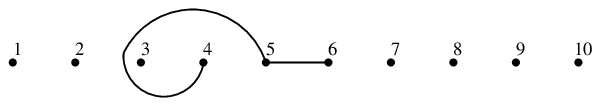}}}$

\medskip
\noindent
$\beta^{-1}_M$ $\hspace{1.7cm}\vcenter{\hbox{\epsfbox{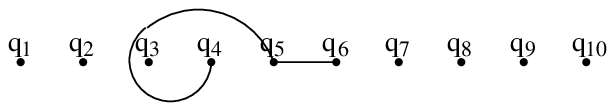}}}$

\medskip
\noindent
$\vp_M(\ell(\ga_{5}))=((Z^2_{45})^{Z^2_{34}} , Z^2_{56})$

\medskip
\noindent
$(\xi_{x'_{11}}) \Psi_{\ga'_{11}}   =  
<4,5>  \De^{\frac{1}{2}}_{I_4I_2} <1>\De<2,3>\De<3,4> \De<5,6>
\De^{\frac{1}{2}}_{I_2\R} <6> \De<4,6> 
\De^2<7,8>\De<9,10>\De<8,9>\De^2<3,4> \beta^{-1}_M = \lcZovere$

\medskip
\noindent
$<4,5>$ $\hspace{.5cm}\vcenter{\hbox{\epsfbox{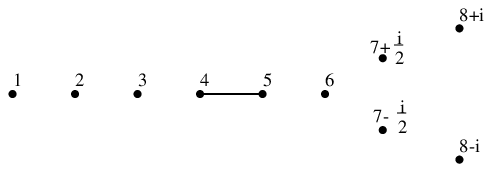}}}$

\medskip
\noindent
$ \De^{\frac{1}{2}}_{I_4I_2} <1>$
$\hspace{.5cm}\vcenter{\hbox{\epsfbox{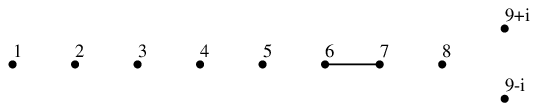}}}$

\medskip
\noindent
$\De<2,3>\De<3,4>$ don't change it

\medskip
\noindent
$\De<5,6>$ $\hspace{.5cm}\vcenter{\hbox{\epsfbox{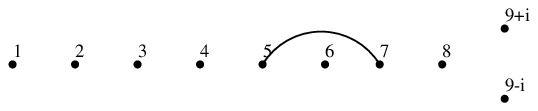}}}$

\medskip
\noindent
$ \De^{\frac{1}{2}}_{I_2\R} <6>$
$\hspace{.5cm}\vcenter{\hbox{\epsfbox{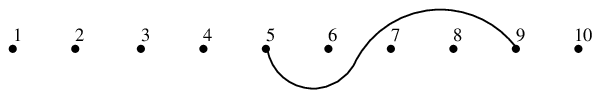}}}$

\medskip
\noindent
$ \De<4,6> $
$\hspace{.5cm}\vcenter{\hbox{\epsfbox{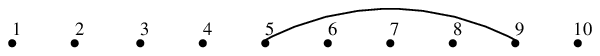}}}$

\medskip
\noindent
$\De^2<7,8>$ doesn't change it

\medskip
\noindent
$\De<9,10>$
$\hspace{.5cm}\vcenter{\hbox{\epsfbox{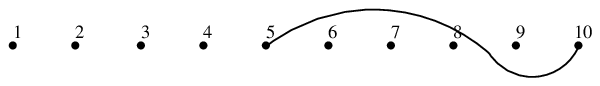}}}$

\medskip
\noindent
$\De<8,9>$
$\hspace{.5cm}\vcenter{\hbox{\epsfbox{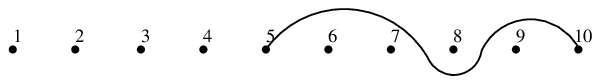}}}$

\medskip
\noindent
$\De^2<3,4>$ doesn't change it

\medskip
\noindent
$\beta^{-1}_M$
$\hspace{.5cm}\vcenter{\hbox{\epsfbox{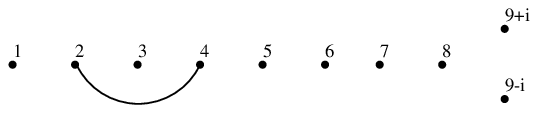}}}$

\medskip
\noindent
$\vp_M(\ell(\ga_{11})) = \Zovere$

\medskip
\noindent
$(\xi_{x'_{13}}) \Psi_{\ga'_{13}}   = <2,3> \De^{\frac{1}{2}}_{I_6I_4} <5> 
\De<4,5>   \De^{\frac{1}{2}}_{I_4I_2} <1> \De<2,3>
\De<3,4>\De<5,6> \De^{\frac{1}{2}}_{I_2\R} <6>\De<4,6>
\De^2<7,8> 
\De<9,10>\De<8,9>\De^2<3,4> \beta^{-1}_M = \stackrel{(4)(6)}{\underline{z}}_{\hspace{-.3cm}3 \; 9}$

\medskip
\noindent
$<2,3>$ $\hspace{.5cm}\vcenter{\hbox{\epsfbox{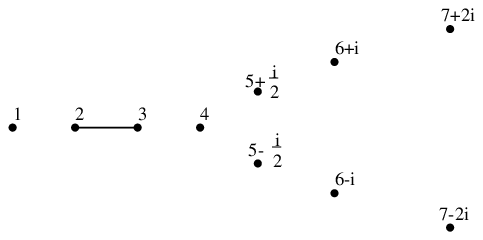}}}$

\medskip
\noindent
$\De^{\frac{1}{2}}_{I_6I_4} <5> $
$\hspace{.5cm}\vcenter{\hbox{\epsfbox{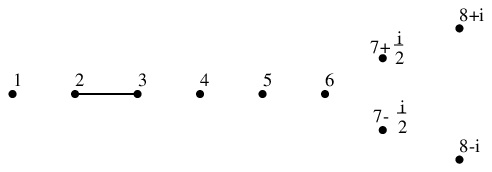}}}$

\medskip
\noindent
$\De<4,5>$ doesn't change it

\medskip
\noindent
$\De^{\frac{1}{2}}_{I_4I_2} <1>$
$\hspace{.5cm}\vcenter{\hbox{\epsfbox{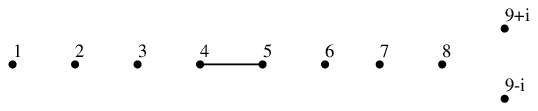}}}$

\medskip
\noindent
$\De<2,3>$ doesn't change it

\medskip
\noindent
$\De<3,4>$
$\hspace{.5cm}\vcenter{\hbox{\epsfbox{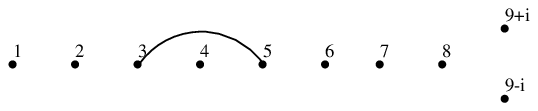}}}$


\medskip
\noindent
$\De<5,6>$
$\hspace{.5cm}\vcenter{\hbox{\epsfbox{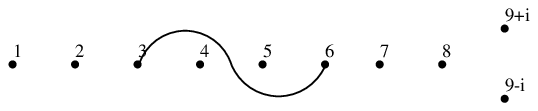}}}$

\medskip
\noindent
$ \De^{\frac{1}{2}}_{I_2\R} <6>$
$\hspace{.5cm}\vcenter{\hbox{\epsfbox{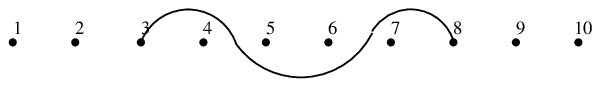}}}$

\medskip
\noindent
$\De<4,6>$
$\hspace{.5cm}\vcenter{\hbox{\epsfbox{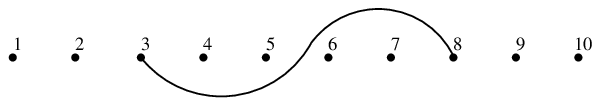}}}$

\medskip
\noindent
$\De^2<7,8>$
$\hspace{.5cm}\vcenter{\hbox{\epsfbox{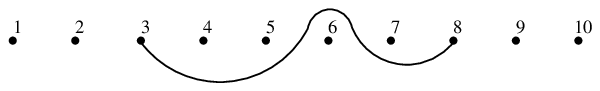}}}$

\medskip
\noindent
$\De<9,10>$ doesn't change it

\medskip
\noindent
$\De<8,9>$
$\hspace{.5cm}\vcenter{\hbox{\epsfbox{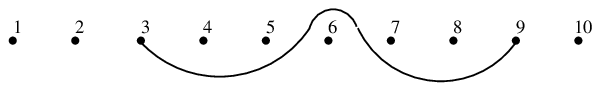}}}$

\medskip
\noindent
$\De^2<3,4>$
$\hspace{.5cm}\vcenter{\hbox{\epsfbox{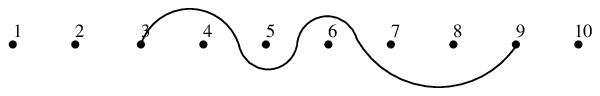}}}$

\medskip
\noindent
$\beta^{-1}_M$
$\hspace{.5cm}\vcenter{\hbox{\epsfbox{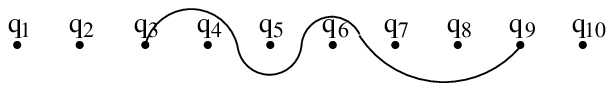}}}$

\medskip
\noindent
$\vp_M(\ell(\ga_{13})) =\stackrel{(4)(6)}{\uZ^2}_{\hspace{-.3cm}3 \; 9}$

Here are the final results of all the other computations:

\medskip
\noindent
$\vp_M(\ell(\ga_{1})) = Z^4_{3 \; 4}$
$\hspace{.5cm}\vcenter{\hbox{\epsfbox{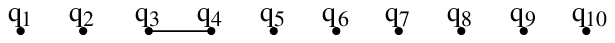}}}$

\medskip
\noindent
$\vp_M(\ell(\ga_{2})) = Z^2_{8 \; 9}$
$\hspace{.5cm}\vcenter{\hbox{\epsfbox{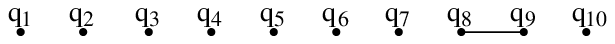}}}$

\medskip
\noindent
$\vp_M(\ell(\ga_{3})) = \bZ^2_{8 \; 10}$
$\hspace{.5cm}\vcenter{\hbox{\epsfbox{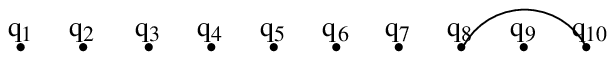}}}$

\medskip
\noindent
$\vp_M(\ell(\ga_{4})) = \uZ^4_{7 \; 9}$
$\hspace{.5cm}\vcenter{\hbox{\epsfbox{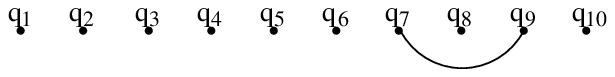}}}$

\medskip
\noindent
$\vp_M(\ell(\ga_{6})) = \left(\bZ_{4 \; 7} \right)^{Z^2_{34}\uZ^2_{79}}$
$\hspace{.5cm}\vcenter{\hbox{\epsfbox{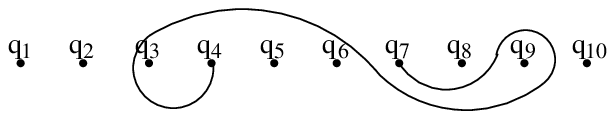}}}$

\medskip
\noindent
$\vp_M(\ell(\ga_{7})) =   \stackrel{(6)}{\uZ^2}_{\hspace{-.1cm}5 \; 9}$
$\hspace{.5cm}\vcenter{\hbox{\epsfbox{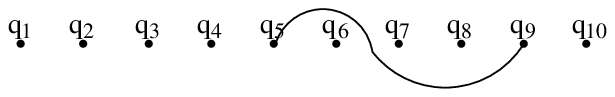}}}$

\medskip
\noindent
$\vp_M(\ell(\ga_{8})) =  \stackrel{(4)}{Z^2}_{\hspace{-.2cm}3 \; 6}$
$\hspace{.5cm}\vcenter{\hbox{\epsfbox{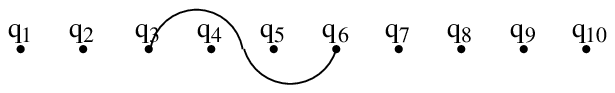}}}$

\medskip
\noindent
$\vp_M(\ell(\ga_{9})) = \uZ^2_{2 \; 6}$
$\hspace{.5cm}\vcenter{\hbox{\epsfbox{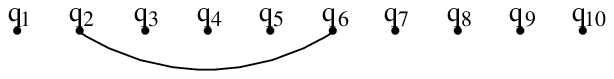}}}$

\medskip
\noindent
$\vp_M(\ell(\ga_{10})) = \uZ_{1 \; 6}$
$\hspace{.5cm}\vcenter{\hbox{\epsfbox{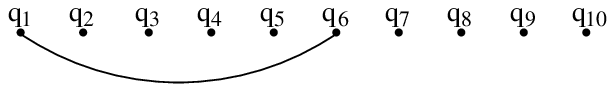}}}$

\medskip
\noindent
$\vp_M(\ell(\ga_{12})) = \bZ_{5 \; 8}$
$\hspace{.5cm}\vcenter{\hbox{\epsfbox{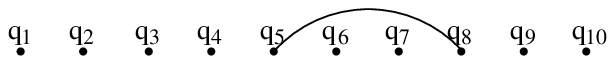}}}$

\medskip
\noindent
$\vp_M (\ell(\ga_{14})) = \stackrel{(6)}{\uZ^2}_{\hspace{-.1cm}2 \; 9} $
$\hspace{.5cm}\vcenter{\hbox{\epsfbox{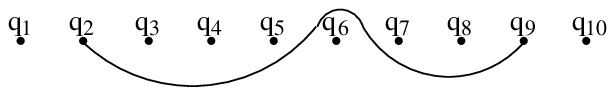}}}$

The sequence of braids that we got is:\\
$Z^4_{3 \; 4} \ , \  Z^2_{8 \; 9} \ , \  \bZ^2_{8 \; 10}    \ , \  
\uZ^4_{7 \; 9} \ , \  \left((Z^2_{45})^{Z^2_{34}}, Z^2_{5 \; 6} \right)
\ , \  \left(\bZ_{4 \; 7} \right)^{Z^2_{34}\uZ^2_{79}} \ , \  
\stackrel{(6)}{\uZ^2}_{\hspace{-.2cm}5 \; 9} \ , \  
\stackrel{(4)}{Z^2}_{\hspace{-.2cm}3 \; 6} \ , \ 
\uZ^2_{2 \; 6} \ , \  \uZ_{1 \; 6} \ , \  
\Zovere
\ , \  \bZ_{5 \; 8} \ , \  \\
\stackrel{(4)(6)}{\uZ^2}_{\hspace{-.3cm}3 \; 9} \ , \ 
\stackrel{(6)}{\uZ^2}_{\hspace{-.2cm}2 \; 9}$
which is the sequence 
that is given by the factorized expression $F_1$. So, as a 
factorized expression, $F_1 = \displaystyle {\prod _{j=1} ^{14}} \varphi _M (\ell (\ga _j))$.

We want to compute $\be _M ^{\vee} (\varphi _M (\ell (\ga_j)))$ for $15 \leq j \leq 28$.
First, we have to compute $\be _{-P} ^{\vee} (\varphi _{-P} (\ell (\tilde \ga_j)))$ 
for $15 \leq j \leq 28$. Then, we have to compute 
$\be _P ^{\vee} (\varphi _P (\ell (\tilde \ga_j T)))$ for $15 \leq j \leq 28$. Finally, we will
compute  
$\be _M ^{\vee} (\varphi _M (\ell ( \ga_j ))) = \be _M ^{\vee} (\varphi _M (\ell ( \tilde \ga_j T T_1 T_2 ))) $
 for $15 \leq j \leq 28$. 

In the computations of $\be _{-P} ^{\vee} (\varphi _{-P} (\ell (\tilde \ga_j)))$, we apply 
$\de _{x_j}$ in a reversed order (see Remark 2.30).
The points $A_{17},A_{23}, A_{19}$ are of the type $a_1$ w.r.t. the point $-P$. 
In the same way as above, we get the following table for the last $14$ points:
\begin{center}
\begin{tabular}{|c|c|c|c|}
\hline
j & $\la _{x_j}$ & $\ep _{x_j}$ & $\de _{x_j}$ \\
\hline
28 & $<7,8>$ & 4 & $\De ^2<7,8>$ \\
27 & $<2,3>$ & 2 & $\De <2,3>$ \\
26 & $<1,2>$ & 2 & $\De <1,2>$ \\
25 & $<3,4>$ & 4 & $\De ^2 <3,4>$ \\
24 & $<5,7>$ & 2 & $\De <5,7>$ \\
23 & $<4,5>$ & 1 & $\De ^{1 \over 2} _{I_2 \R} <4>$ \\
22 & $<5,6>$ & 2 & $\De <5,6>$ \\
21 & $<3,4>$ & 2 & $\De <3,4>$ \\
20 & $<6,7>$ & 2 & $\De <6,7>$\\
19 & $<7,8>$ & 1 & $\De ^{1\over 2} _{I_4 I_2} <7>$\\
18 & $<2,3>$ & 2 & $\De <2,3>$\\
17 & $<1,2>$ & 1 & $\De ^{1 \over 2} _{I_6 I_4} <1>$\\
16 & $<2,3>$ & 2 & $\De <2,3>$\\
15 & $<3,4>$ & 2 & $\De <3,4>$\\
\hline
\end{tabular}
\end{center}

We compute $\be _{-P} ^{\vee} (\varphi _{-P} (\ell (\tilde \ga_j)))$ for $j=28,\cdots,15$.
We give here just the computations for: $j=24,18,16$.
\def\tg{\tilde{\gamma}}

\noindent
L.V.C. $(\tg_{24})\beta_{-P} = <5,7>\De^2<3,4>\De<1,2>\De<2,3>\De^2<7,8> = 
(\underline{z}_{5 \; 6}, (z_{6 \; 7})^{Z^2_{7 \; 8}})$

\medskip
\noindent
$<5,7>$
$\hspace{.5cm}\vcenter{\hbox{\epsfbox{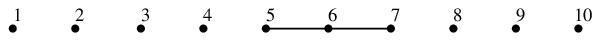}}}$

\medskip
\noindent
$ \De^2<3,4>\De<1,2>\De<2,3>$don't change it

\medskip
\noindent
$\De^2<7,8>$
$\hspace{.5cm}\vcenter{\hbox{\epsfbox{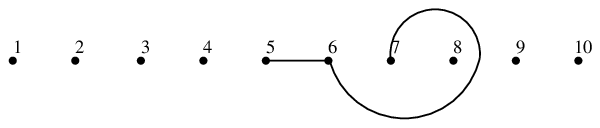}}}$

\medskip
\noindent
$\beta^\vee_{-P} (\vp_{-P}(\ell(\tg_{24}))) = (\uZ^2_{5 \; 6}, (Z^2_{6 \; 7})^{Z^2_{7 \; 8}})$

\medskip
\noindent
L.V.C.$ (\tg_{18})\beta_{-P} = <2,3>  \De^{\frac{1}{2}}_{I_4I_2} <7> \De<6,7>\De<3,4>\De<5,6> \De^{\frac{1}{2}}_{I_2\R} <4> \De<5,7>\De^2<3,4>\De<1,2>\De<2,3>
\De^2<7,8> = \stackrel{(3)}{\underline{z}}_{1 \ 6}$

\medskip
\noindent
$<2,3>$
$\hspace{.5cm}\vcenter{\hbox{\epsfbox{p733.ps}}}$

\medskip
\noindent
$\De^{\frac{1}{2}}_{I_4I_2} <7>$
$\hspace{.5cm}\vcenter{\hbox{\epsfbox{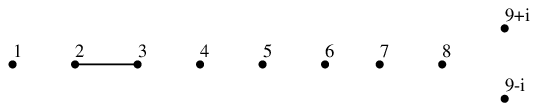}}}$

\medskip
\noindent
$ \De<6,7>$ doesn't change it\\

\bigskip
\noindent
$\De<3,4>$
$\hspace{.5cm}\vcenter{\hbox{\epsfbox{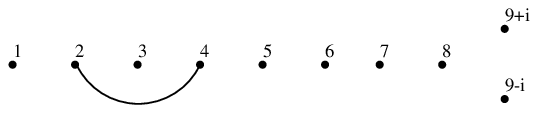}}}$

\medskip
\noindent
$ \De<5,6>$ doesn't change it

\medskip
\noindent
$ \De^{\frac{1}{2}}_{I_2\R} <4>$
$\hspace{.5cm}\vcenter{\hbox{\epsfbox{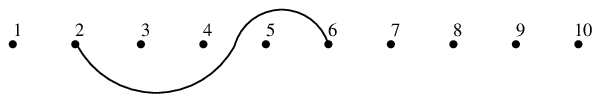}}}$

\medskip
\noindent
$  \De<5,7>$
$\hspace{.5cm}\vcenter{\hbox{\epsfbox{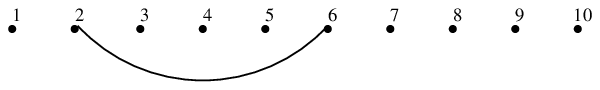}}}$

\medskip
\noindent
$ \De^2<3,4>$ doesn't change it

\medskip
\noindent
$\De<1,2>$
$\hspace{.5cm}\vcenter{\hbox{\epsfbox{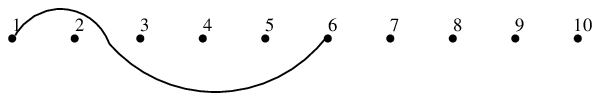}}}$

\medskip
\noindent
$\De<2,3>$
$\hspace{.5cm}\vcenter{\hbox{\epsfbox{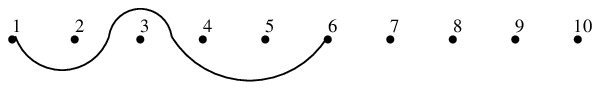}}}$

\noindent
$\De^2<7,8>$ doesn't change it

\medskip
\noindent
$\beta^\vee_{-P} (\vp_{-P}(\ell(\tg_{18}))) = \stackrel{(3)}{\uZ^2}_{
\vspace{-.2cm}1 \ 6}$\\

\medskip
\noindent
L.V.C.$ (\tg_{16})\beta_{-P} = <2,3>  \De^{\frac{1}{2}}_{I_6I_4} <1> 
\De<2,3> \De^{\frac{1}{2}}_{I_4I_2} <7> \De<6,7>\De<3,4>\De<5,6> \De^{\frac{1}{2}}_{I_2\R} <4> \De<5,7>\De^2<3,4>\De<1,2>\De<2,3>
\De^2<7,8> = \lcZoverfs$

\medskip
\noindent
$<2,3>$
$\hspace{.5cm}\vcenter{\hbox{\epsfbox{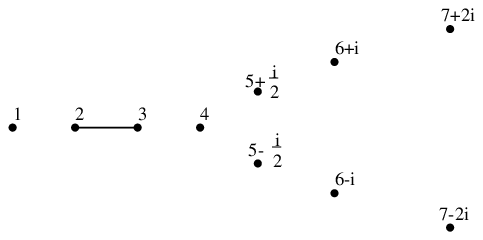}}}$

\medskip
\noindent
$ \De^{\frac{1}{2}}_{I_6I_4} <1> $
$\hspace{.5cm}\vcenter{\hbox{\epsfbox{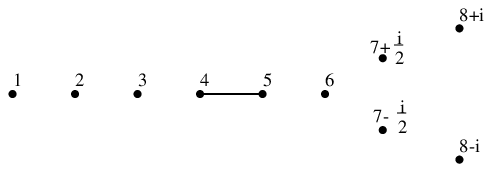}}}$

\medskip
\noindent
$ \De<2,3>$ doesn't change it

\medskip
\noindent
$ \De^{\frac{1}{2}}_{I_4I_2} <7> $
$\hspace{.5cm}\vcenter{\hbox{\epsfbox{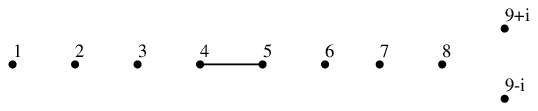}}}$

\medskip
\noindent
$\De<6,7>$ doesn't change it

\medskip
\noindent
$\De<3,4>$
$\hspace{.5cm}\vcenter{\hbox{\epsfbox{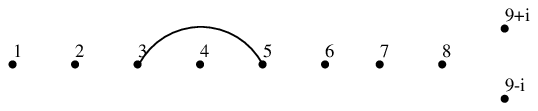}}}$

\medskip
\noindent
$\De<5,6>$
$\hspace{.5cm}\vcenter{\hbox{\epsfbox{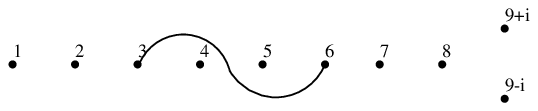}}}$

\medskip
\noindent
$\De^{\frac{1}{2}}_{I_2\R} <4>$
$\hspace{.5cm}\vcenter{\hbox{\epsfbox{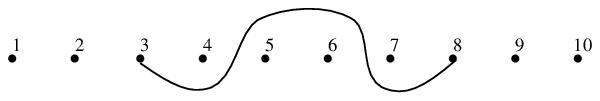}}}$

\medskip
\noindent
$\De<5,7>$
$\hspace{.5cm}\vcenter{\hbox{\epsfbox{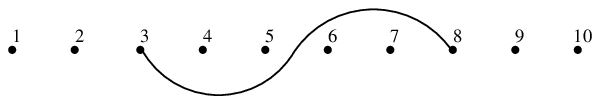}}}$

\medskip
\noindent
$\De^2<3,4>$
$\hspace{.5cm}\vcenter{\hbox{\epsfbox{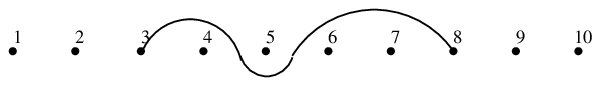}}}$

\medskip
\noindent
$\De<1,2>$ doesn't change it

\medskip
\noindent
$\De<2,3>$
$\hspace{.5cm}\vcenter{\hbox{\epsfbox{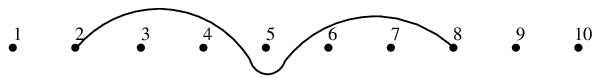}}}$

\medskip
\noindent
$\De^2<7,8>$
$\hspace{.5cm}\vcenter{\hbox{\epsfbox{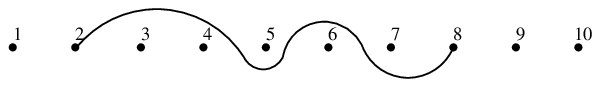}}}$

\medskip
\noindent
$\beta^\vee_{-P} (\vp_{-P}(\ell(\tg_{16}))) =\Zoverfs$

Here are the final results of all the other computations:\\


\medskip
\noindent
$\vcenter{\hbox{\epsfbox{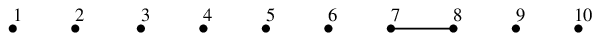}}}$\hspace{5cm}
$\beta^\vee_{-P} (\vp_{-P}(\ell(\tg_{28}))) = Z^4_{7 \; 8}$

\medskip
\noindent
$\vcenter{\hbox{\epsfbox{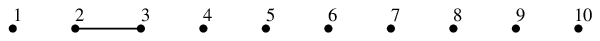}}}$\hspace{5cm}
$\beta^\vee_{-P} (\vp_{-P}(\ell(\tg_{27}))) = Z^2_{2 \; 3}$

\medskip
\noindent
$\vcenter{\hbox{\epsfbox{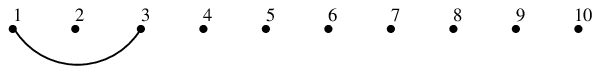}}}$\hspace{5cm}
$\beta^\vee_{-P} (\vp_{-P}(\ell(\tg_{26}))) = \uZ^2_{1 \; 3}$

\medskip
\noindent
$\vcenter{\hbox{\epsfbox{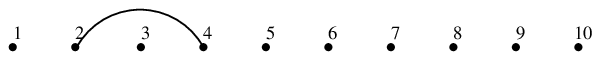}}}$\hspace{5cm}
$\beta^\vee_{-P} (\vp_{-P}(\ell(\tg_{25}))) = \bZ^4_{2 \; 4}$

\medskip
\noindent
$\vcenter{\hbox{\epsfbox{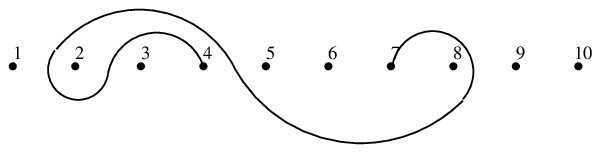}}}$\hspace{5cm}
$\beta^\vee_{-P} (\vp_{-P}(\ell(\tg_{23}))) = 
\left(\uZ_{4 \; 7} \right)^{Z^2_{78}\bZ^2_{24}}$

\medskip
\noindent
$\vcenter{\hbox{\epsfbox{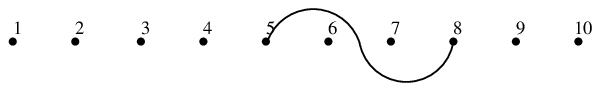}}}$\hspace{5cm}
$\beta^\vee_{-P} (\vp_{-P}(\ell(\tg_{22}))) = 
\stackrel{(6)}{Z^2}_{\hspace{-.2cm}5 \; 8}$

\medskip
\noindent
$\vcenter{\hbox{\epsfbox{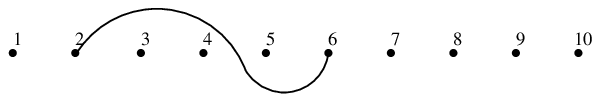}}}$\hspace{5cm}
$\beta^\vee_{-P} (\vp_{-P}(\ell(\tg_{21}))) = 
\Zoverfivetwosix$

\medskip
\noindent
$\vcenter{\hbox{\epsfbox{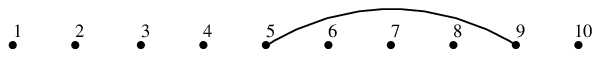}}}$\hspace{5cm}
$\beta^\vee_{-P} (\vp_{-P}(\ell(\tg_{20}))) = 
\bZ^2_{5 \; 9}$

\medskip
\noindent
$\vcenter{\hbox{\epsfbox{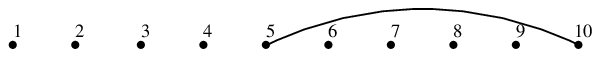}}}$\hspace{5cm}
$\beta^\vee_{-P} (\vp_{-P}(\ell(\tg_{19}))) = 
\bZ_{5 \; 10}$

\medskip
\noindent
$\vcenter{\hbox{\epsfbox{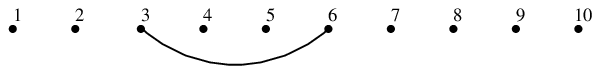}}}$\hspace{5cm}
$\beta^\vee_{-P} (\vp_{-P}(\ell(\tg_{17}))) = 
\uZ_{3 \; 6}$

\medskip
\noindent
$\vcenter{\hbox{\epsfbox{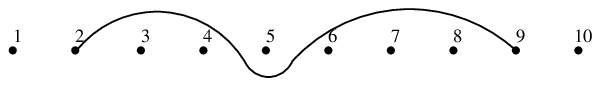}}}$\hspace{5cm}
$\beta^\vee_{-P} (\vp_{-P}(\ell(\tg_{15}))) = 
\Zoverfive$
\\

Now, we give the results of ${\rm L.V.C.} (\tilde \ga _j T) \be _P$ for 
$ j= 28, \cdots, 15 \ \ $ from [Am].
$${\rm L.V.C.} (\tilde \ga _j T) \be _P ={\rm L.V.C.} (\tilde \ga _j)\Psi _T \be _P = {\rm L.V.C.} (\tilde \ga _j) \be _{-P} (\be _{-P} ^{-1}\Psi _T \be _P), \ \ j=28, \cdots, 15.$$
Because $S_2$ does not have any singularities at $\infty$, 
$\be _{-P} ^{-1}\Psi _T \be _P$ is a $180 ^{\circ}$ rotation centered   at $5 {1\over 2}$.

So we apply $180 ^{\circ}$ rotation on  ${\rm L.V.C.} (\tilde \ga _j) \be _{-P}$ to get 
${\rm L.V.C.} (\tilde \ga _j T) \be _P$:

\medskip
\noindent
L.V.C.$(\tg_{28}T)\beta_P = 
\vcenter{\hbox{\epsfbox{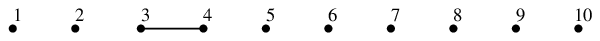}}}
  = z_{3 \; 4}$ 

\medskip
\noindent L.V.C.$(\tg_{27}T)\beta_P = \vcenter{\hbox{\epsfbox{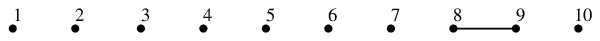}}}
 = z_{8 \; 9}$ 

\medskip
\noindent L.V.C.$(\tg_{26}T)\beta_P = 
\vcenter{\hbox{\epsfbox{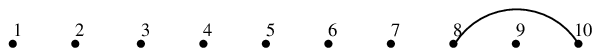}}}
 = \bar{z}_{8 \; 10}$

\medskip
\noindent L.V.C.$(\tg_{25}T)\beta_P = 
\vcenter{\hbox{\epsfbox{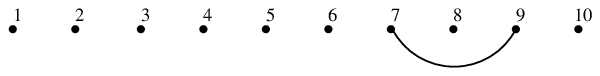}}}
= \underline{z}_{7 \; 9}$ 

\medskip
\noindent L.V.C.$(\tg_{24}T)\beta_P = 
\vcenter{\hbox{\epsfbox{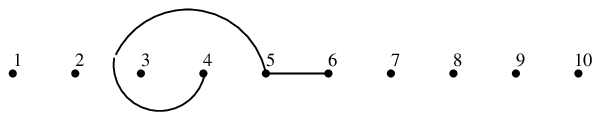}}}
 = ((z_{4\;5})^{Z^2_{3\;4}} , 
z_{5 \; 6})$

\medskip
\noindent L.V.C.$(\tg_{23}T)\beta_P = 
\vcenter{\hbox{\epsfbox{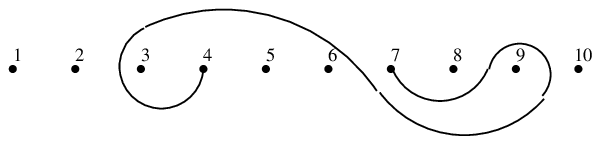}}}
 = (\bar{z}_{4\;7})^{Z^2_{3\;4}  
\uZ^2_{7 \; 9}}$

\medskip
\noindent L.V.C.$(\tg_{22}T)\beta_P = 
\vcenter{\hbox{\epsfbox{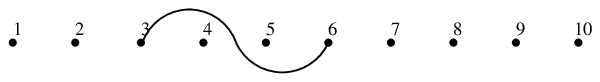}}}
 = 
\stackrel{(4)}{z}_{3 \; 6}$

\medskip
\noindent L.V.C.$(\tg_{21}T)\beta_P = 
\vcenter{\hbox{\epsfbox{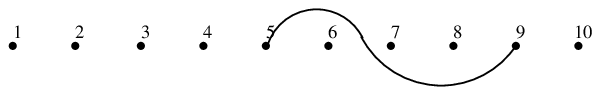}}}
 = 
\stackrel{(6)}{\underline{z}}_{5 \; 9}$

\medskip
\noindent L.V.C.$(\tg_{20}T)\beta_P = 
\vcenter{\hbox{\epsfbox{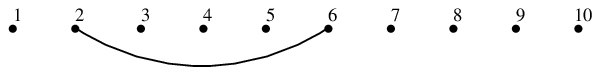}}}
 = 
\underline{z}_{2 \; 6}$

\medskip
\noindent L.V.C.$(\tg_{19}T)\beta_P = 
\vcenter{\hbox{\epsfbox{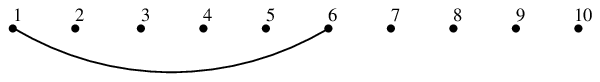}}}
 = 
\underline{z}_{1 \; 6}$

\medskip
\noindent L.V.C.$(\tg_{18}T)\beta_P = 
\vcenter{\hbox{\epsfbox{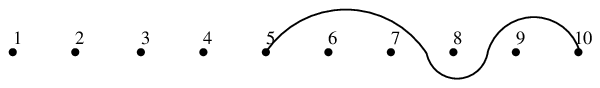}}}
 = 
\zovereight$

\medskip
\noindent L.V.C.$(\tg_{17}T)\beta_P = 
\vcenter{\hbox{\epsfbox{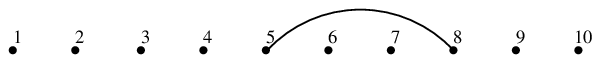}}}
 = 
\bar{z}_{5 \; 8}$

\medskip
\noindent L.V.C.$(\tg_{16}T)\beta_P =
\vcenter{\hbox{\epsfbox{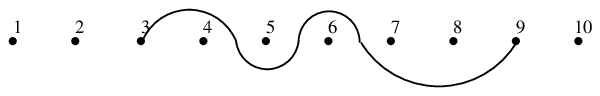}}}
 = 
\stackrel{(4)(6)}{\underline{z}}_{3 \; 9}$

\medskip
\noindent L.V.C.$(\tg_{15}T)\beta_P = 
\vcenter{\hbox{\epsfbox{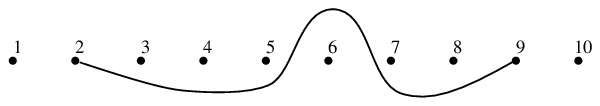}}}
= 
\stackrel{(6)}{\underline{z}}_{2 \; 9}$

\def\vp{\varphi}
Moreover, $\be _P ^{\vee} (\varphi _P (\ell (\tilde \ga _j T))) = \De ^{\ep _{x_j}} < {\rm L.V.C.} (\tilde \ga _j T) \be _P>\ \ $ for $j =28, \cdots, 15$. So:
\\
\\
$\beta^\vee_P(\vp_P(\ell(\tg_{28}T))) = Z^4_{3\;4}$\\
\\
$\beta^\vee_P(\vp_P(\ell(\tg_{27}T))) =Z^2_{8 \; 9}$\\
\\
$\beta^\vee_P(\vp_P(\ell(\tg_{26}T))) = \bZ^2_{8 \; 10} $\\
\\
$\beta^\vee_P(\vp_P(\ell(\tg_{25}T))) = \uZ^4_{7 \; 9} $\\
\\
$\beta^\vee_P(\vp_P(\ell(\tg_{24}T))) =((Z^2_{4\;5})^{Z^2_{34}} , Z^2_{5 \; 6})$\\ 
\\
$\beta^\vee_P(\vp_P(\ell(\tg_{23}T))) =\left(\bZ_{4 \; 7} \right)^{Z^2_{34}\uZ^2_{79}}$\\
\\
$\beta^\vee_P(\vp_P(\ell(\tg_{22}T))) = \stackrel{(4)}{Z^2}_{\hspace{-.2cm}3 \; 6}$ \\
\\
$\beta^\vee_P(\vp_P(\ell(\tg_{21}T))) =  \stackrel{(6)}{\uZ^2}_{\hspace{-.2cm}5 \; 9}$\\
\\
$\beta^\vee_P(\vp_P(\ell(\tg_{20}T))) = \uZ^2_{2 \; 6}$\\
\\
$\beta^\vee_P(\vp_P(\ell(\tg_{19}T))) = \uZ_{1 \; 6}$\\
\\
$\beta^\vee_P(\vp_P(\ell(\tg_{18}T))) = \Zovereight$\\
\\
$\beta^\vee_P(\vp_P(\ell(\tg_{17}T))) = \bZ_{5 \; 8}$\\
\\
$\beta^\vee_P(\vp_P(\ell(\tg_{16}T))) = \stackrel{(4)(6)}{\uZ^2}_{\hspace{-.3cm}3 \; 9}$\\
\\
$\beta^\vee_P(\vp_P(\ell(\tg_{15}T))) = \stackrel{(6)}{\uZ^2}_{\hspace{-.2cm}2 \; 9}$

By Figure \ref{s2fig}, L-pair$(x(l_2 \cap l_3)) = (2,3)$, 
L-pair$(x(l_6 \cap l_7)) = (9,10)$. These two points are of the type $c$. 
By Theorem 2.24, $ \be _P ^{-1} \Psi _{T_1 T_2} \be _M = (\De <2,3> \De <9,10> )^{-1}$. 

Let us denote $\De <2,3> \De <9,10>$ by $\rho$. 

Now, for $j = 28, \cdots, 15\ \ $ we have the following:\\
$\vp_M(\ell(\tg_jTT_1T_2)) = \Delta^{\epsilon_{x_j}}<{\rm L.V.C.} (\tg_jTT_1T_2)> = \Delta^{\epsilon_{x_j}}< {\rm L.V.C.} (\tg_jT)\Psi_{T_1T_2}> = \\
\Delta^{\epsilon_{x_j}}< {\rm L.V.C.} (\tg_jT)\beta_P (\beta^{-1}_P \Psi_{T_1T_2}\beta_M)\beta^{-1}_M> = \Delta^{\epsilon_{x_j}}< {\rm L.V.C.} (\tg_jT)\beta_P \rho^{-1}\beta^{-1}_M> =\\
(\beta^{-1}_M)^\vee(\rho^{-1})^\vee \Delta^{\epsilon_{x_j}}< {\rm L.V.C.} (\tg_j T)
\beta_P> = (\beta^{-1}_M)^\vee(\rho^{-1})^\vee (\beta^\vee_P 
(\vp_P(\ell(\tg_j T))))$.
\\
So by comparing the beginning and the end we have:
$$\varphi _M (\ell (\ga _j)) = 
(\be _M ^{-1})^{\vee} \be _P ^{\vee}
 ((\varphi _P (\ell (\tilde \ga _j T)))^{\rho^{-1}} )\ , 
\qquad  j= 28, \cdots, 15.$$ 
Since $\be _M (q_i)=i$, one can easily see that as a factorized expression:
$$\displaystyle F_1 = {\prod _{j=28} ^{15}} (\be _M ^{-1})^{\vee} \be _P ^{\vee} (\varphi _P (\ell (\tilde \ga _j T))).$$  

So, as a factorized expression:
$$\displaystyle F_1^{\rho^{-1}} = {\prod _{j=28} ^{15}} (\be _M ^{-1})^{\vee} \be _P ^{\vee} ((\varphi _P (\ell (\tilde \ga _j T)))^{\rho^{ -1}} ) = {\prod _{j=28} ^{15}} \varphi _M (\ell (\ga _j)).$$  
Hence, the braid monodromy w.r.t. $E \times D$ 
is $\varphi _M =F_1 \cdot F_1 ^{\rho^{-1}}$.
\hfill
$\Box$

\begin{\bib}{10}
\bibitem[A] {} Artin, E., {\it Theory of braids}, Ann. of
Math. {\bf 48} (1947), 101-126.
\bibitem[Am]{Am} Amram, M., {\it Braid group and braid monodromy}, M.Sc. thesis (1995)
[in Hebrew].
\bibitem[Ha]{Ha} { Hartshorne, R., {\it Algebraic Geometry}, Graduate Texts in 
   Math., Vol. {\bf 52} (1977).}
\bibitem[K]{K} Kendig, K., {\it Elementary Algebraic Geometry}, Graduate Texts in 
   Math., Vol. {\bf 44} (1977).
\bibitem[MoTe1]{MoTe1} Moishezon, B.,  Teicher, M., {\it Braid group 
   technique in complex geometry I, Line arrangements in $\C\PP ^2$}, 
   Contemporary Math. {\bf 78}, 425-555 (1988).  
\bibitem[MoTe2]{MoTe2} Moishezon, B., Teicher, M., {\it Braid group 
   technique in complex geometry II, From arrangements of lines and conics 
   to cuspidal curves}, Algebraic Geometry, Lect. Notes in Math. {\bf 1479} 
   (1990).
\bibitem[MoTe3]{MoTe3} Moishezon, B., Teicher, M., {\it Braid group 
   technique in complex geometry III, Projective degeneration of $V_3$}, 
   Contemporary Math. {\bf 162}, 313-332 (1994). 
\bibitem[MoTe4]{MoTe4} Moishezon, B., Teicher, M., {\it Braid group 
   technique in complex geometry IV, Braid monodromy of the branch curve 
   $S_3$ of $V_3 \to \C \PP ^2$ and application to $\pi _1 (\C\PP ^2 -S_3, *)$}, 
   Contemporary Math. {\bf 162}, 332-358 (1994). 
\bibitem[VK]{VK} Van Kampen, E.R., {\it On the fundamental group of 
   an algebraic curve}, Amer. J. Math. {\bf 55}, 255-260 (1933). 
\bibitem[Zu]{Zu} Zur-Cycowicz, S., {\it Braid monodromy}, 
   M.Sc. thesis (1993) [in Hebrew].
\end{\bib}

\end{document}